\documentclass[journal]{IEEEtran}
\usepackage{cite}

\usepackage{graphicx}

\usepackage{amssymb}
\usepackage[cmex10]{amsmath}

\usepackage{caption}
\usepackage[font=footnotesize]{subfig}

\usepackage{diagbox}
\usepackage[flushleft]{threeparttable}
\usepackage{tabularx}
\usepackage{multirow}
\usepackage{soul}
\usepackage{color}
\usepackage{optidef}
\usepackage{enumitem}

\begin{document}
%
\title{Deriving DERs VAR-Capability Curve at TSO-DSO Interface to Provide Grid Services}


\author{A.~Singhal,~\IEEEmembership{Member,~IEEE,}
A. K.~Bharati,~\IEEEmembership{Member,~IEEE,}
        V.~Ajjarapu,~\IEEEmembership{Fellow,~IEEE}%
}%


%


\maketitle

\begin{abstract}
The multitudes of inverter-based distributed energy resources (DERs) can be envisioned as distributed reactive power (var) devices (\textit{mini-SVCs}) that can offer var flexibility at TSO-DSO interface. To facilitate this vision, a systematic methodology is proposed to derive an aggregated var capability curve of a distribution system with DERs at the substation level, analogous to a conventional bulk generator. Since such capability curve will be contingent to the operating conditions and network constraints, an optimal power flow (OPF) based approach is proposed that takes inverter headroom flexibility, unbalanced nature of system and coupling with grid side voltage into account along with changing operating conditions. Further, the influence of several factors such as compliance to IEEE 1547 on the capability curve is thoroughly investigated on an IEEE 37 bus and 123 bus distribution test system along with unbalanced DER proliferation. Validation with nonlinear analysis is presented along with demonstration of a scenario with T-D co-simulation.


\end{abstract}

\IEEEpeerreviewmaketitle

\section{Introduction}
\IEEEPARstart{R}{eactive} power (var) balance plays a vital role in maintaining transmission grid resiliency and, availability of sufficient var capability is often considered an indicator of voltage security \cite{bao_online_2003,song_reactive_2003}. The var related ancillary services have been mainly achieved by large synchronous generators and other strategically deployed var devices such as static synchronous compensator (STATCOM) and static var compensator (SVC). However, a growing footprint of distributed energy resources (DERs) is replacing fossil fuel based generation that may result in shortage of regional var availability \cite{goergens_determination_2015,barth_technical_2013}. It has initiated a discussion on utilizing DERs as alternative sources in the future grid, along with bulk generation plants, to provide essential ancillary services to the grid such as ramping requirements, ensuring adequate inertia, and maintaining var balance \cite{keane_state---art_2013,perez-arriaga_transmission_2016}. This paper concerns to the latter topic with DER as focus. As voltage insecurity or voltage instability is usually load driven and DERs are closer to the load centers compared to the generating stations, DERs may serve some part of the var requirement.

Most of the inverter-based DERs have capability to independently control the real power and var. Much of the extant literature focuses on utilizing the var control potential of DERs to improve the performance of distribution systems (DS) i.e. voltage challenges \cite{singhal_real-time_2018,zhu_fast_2016}, loss minimization \cite{zhang_optimal_2015} and other such indices. 
{\color{blue}However, utilization of DERs' var potential for the benefit of the transmission systems (TS) is being explored more recently.

We present a methodology that enables thousands of DER devices with var control capability can be utilized as the geographically distributed var resources (\textit{`mini- SVCs'}). }These (\textit{`mini- SVCs'}) can provide enhanced flexibility options to the transmission system operators (TSOs) either directly if located close to the substation or indirectly by reducing the demand by the distribution system operators (DSO), if coordinated properly. 
{\color{blue}Our assertion} is founded on the following reasoning: 1) The inverter-based DERs can inject/absorb var via fast local volt/var controls \cite{singhal_real-time_2018,zhu_fast_2016}, thus can provide a significant amount of fast and continuous capacitive/inductive var support, locally or at an aggregated level; 2) The  proposition of DERs' var provision is gaining strength with revised DER integration standards such as IEEE1547-2018 \cite{noauthor_ieee_2018}, California Rule 21, Hawaii Rule 14 \cite{noauthor_impact_2018} and Germany grid codes \cite{barth_technical_2013} that have made it obligatory for DERs to provide var support for grid requirements; 3) The local and distributed nature of DER makes it a suitable contestant for var provision. In fact, an assessment study for East Denmark identifies DERs var provision scheme economically competitive to conventional dynamic var devices. \cite{barth_technical_2013}; 4) The required infrastructure and protocol for DSO-TSO interaction has started gaining attention e.g. 
some TSOs in Europe 
are implementing a payment structure for voltage control where DSOs can participate in var provision based on the day-ahead reactive power plans sent out by TSO \cite{marten_analysis_2013}. 

Thus, in this new environment of DERs, a consensus emerge from literature that motivates TSOs to consider DERs var flexibility in their optimization.
Previous works have aggregated the capability of asynchronous generators or DFIG for large wind farms without considering DS constraints\cite{konopinski_extended_2009,cuffe_transmission_2012}
 Reference \cite{kundu_approximating_2018} has attempted to approximate the DER flexibility using geometric approach without physical network constraints. \cite{marten_optimizing_2014,goergens_determination_2015,kaempf_reactive_2014} introduce the optimization based approach with focus on the TSO-DSO interaction. 
{\color{blue}There are some papers that discuss estimating the DER capability to provide ancillary services. Reference \cite{8442917} has discussed a mixed integer linear program (MILP) at the transmission system level that is designed for demand reduction coupled with OLTC operation with an OLTC headroom to account for distribution system flexibility, however, does not model the distribution systems. Reference \cite{8291006} has discussed an interval constrained {\color{blue}power flow} optimization technique that estimates the active power and reactive power flexibility at the boundary of the TSO-DSO. Reference \cite{9543347} has discussed a method to estimate the flexibility of the distribution system, however the high-level framework is similar to Reference \cite{8291006} but the optimization formulation of the problem is different and uses the PV and energy storage DERs to achieve the DER flexibility. Reference \cite{CAPITANESCU2018226} provides a method of estimating the maximum and minimum real and reactive power flexibility from the distribution system with DERs. However, these methods do not consider inverter and network constraints which may become critical and binding when we aggregate large number of devices. References \cite{9295337} have discussed methods to to account for network and inverter constraints in flexibility estimation.
Similarly, \cite{Kara} also discusses the var support by distribution systems to the bulk grid with network constraints.

However, we identify following main limitations of the existing work: (1) These works do not account and analyze the impact of unbalance in distribution systems on the flexibility estimation, (2) The flexibility is not characterized as function of inverter headroom as per the IEEE 1547 standard and, (3) these papers do not discuss how the capability curves could be impacted by various practical aspects of distribution system such as grid side voltage, unbalanced DER distribution.} The present paper is an extension of the conference paper, \cite{singhal_framework_2018}, where the extension discusses the implications of considering headroom while operating the DERs and its impact on the TSO's ability to utilize the DSO's var availability. A more comprehensive, detailed understanding of the DER var estimation at the TSO-DSO interface along with validation is presented in this extended paper. This paper also discusses the potential impact of these flexiblity curves on transmission grid.
\subsection{IEEE 1547 and Headroom Correlation}
The IEEE 1547-2018 standard \cite{noauthor_ieee_2018} mandates at least 44 \% of the inverter rating to be available for reactive power modulation. For a solar PV DER, the implication is that for 90\% of the time, when the solar insolation is not peaking, the inverter will be oversized by ~ 11\% that means higher investment. Therefore, it is expected that most of the inverters will not be oversized and will need ~10\% of real power headroom during the peak solar PV generation (minimum headroom) to  allow for 44\% of the inverter rating for reactive power modulation.

\begin{figure}
	\centering
    \vspace{-4mm}
	\includegraphics[trim=0.5in 2in 1.5in 2.5in,width=3.4in]{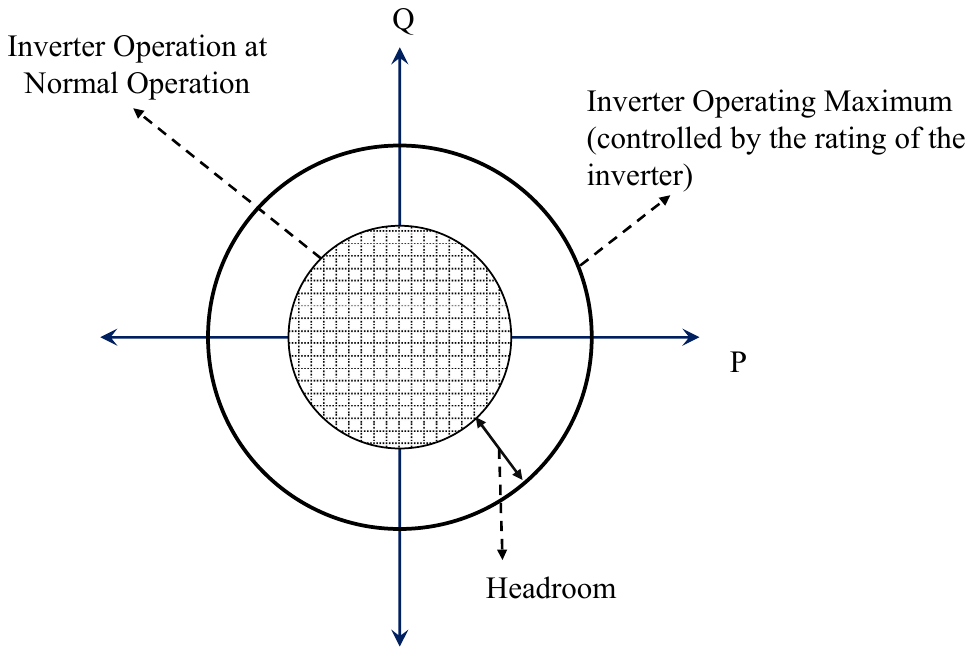}
    \vspace{-3mm}
    \caption {Difference in Inverter Operation and Maximum Rating that Creates a Headroom}
    \label{fig:headroom}
\end{figure}

The capacity of the inverters cannot be utilized if there is no headroom. {\color{blue}Another implication of the var capability is not operating the inverter in unity power factor mode, that results in poor power factor and its adverse effects.} Therefore, IEEE 1547 mandates a 44\% reactive power modulation capability under all operating conditions.  \figurename \ref{fig:headroom} demonstrates how headroom is formed and if it is maintained then the inverter can supply additional real and reactive power within its rating under emergency conditions.

{\color{blue}
\subsection{Key Contributions}
The present work provides following unique contribution in DER capability/flexibility aggregation: 1) The proposed methodology provides an aggregated net \textit{Q-capability} curve as function of real power headroom (complying to the IEEE 1547-2018 standard) resembling a virtual conventional generator capability curve. This will enable TSO to model both P and Q flexibility as resources from DS into their planning.
2) The proposed approach captures the unbalanced nature of distribution system which is shown to be a crucial factor affecting flexibility. (3) To provide useful and comprehensive insight, the influence of several factors on aggregated capability is investigated such as grid side voltage, daily load profile, revised integration standard 1547, inverter sizing etc. 4) Few simple applications of the derived DER flexibility for the transmission grid has been demonstrated using T-D co-simulation. T-D cosimulation allows to observe the impact of flexibility on TS while ensuring that DS operational constraints are not violated.}


Section II sets up the conceptual framework for var provision. Section III builds capability curve characterization for which an OPF based process is outlined in Section IV. Section V presents DS case studies with discussion on the IEEE1547 standard compliance and unbalanced DER proliferation in DS. Section VI demonstrates the impact of DER var provision on the grid via T-D cosimulation.

\section{Overall Conceptual VAR Support Framework }
\begin{figure}
	\centering
    \vspace{-4mm}
	\includegraphics[trim=-0.5in 0in 0in 0in,width=3.4in]{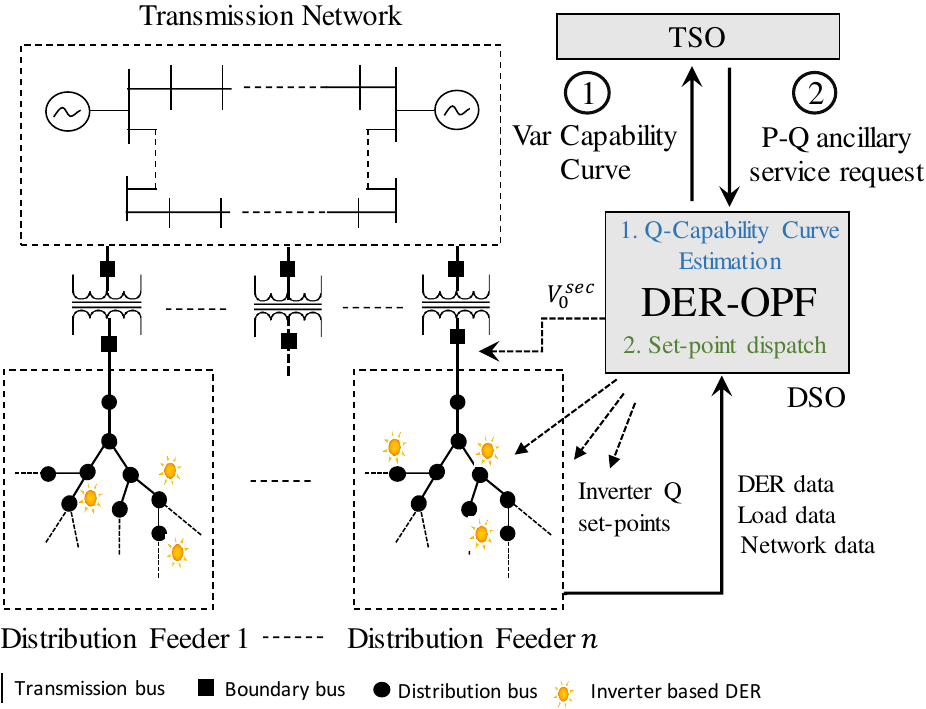}
    \vspace{-3mm}
    \caption {The proposed DER var support framework for an integrated T-D system which has two major functionalities for DSO. Functionality 1, providing 'Var Capability Curve' to TSO is the focus of this paper. }
    \label{fig:overall}
\end{figure}

 \figurename \ref{fig:overall} depicts the overall framework of providing DERs' var support to the grid in an integrated T-D system, proposed in this work. Consider a transmission grid which is connected to multiple DS with high penetration of inverter-based DERs. In this study, distributed solar photovoltaic (PV) are considered as DERs. The whole physical system can be seen in three parts i.e. transmission grid, boundary buses (substation) and the distribution buses with DERs. In this framework, we envision an aggregator entity DSO at substation level which exchange information with the TSO and the DER devices. 
As shown in the \figurename \ref{fig:overall}, the framework consists of two major functions performed by the DSO. However, in this paper we only focus on the first function that is to dynamically aggregate the \textit{net var capability curve} of the DS at the substation level in every ~10-15 minutes time scale based on short-term forecast and send it to the TSO to include it in their planning and operational activities. Here we assume that the TSO has its own planning and control methods to request var support from the DSO in case of emergency. The second function of DSO is to dispatch optimal inverter var set-points to individual DER devices in order to meet the var support requested by the grid, however, in this work we do not provide details of this functionality. The scope of this paper is to focus on developing a general framework to aggregate DER var capability. 
 

\section{Capability Curve Characterization }
\begin{figure}
\centering
\vspace{-4mm}
\includegraphics[trim=-0in 0in 0in 0in,width=3.5in]{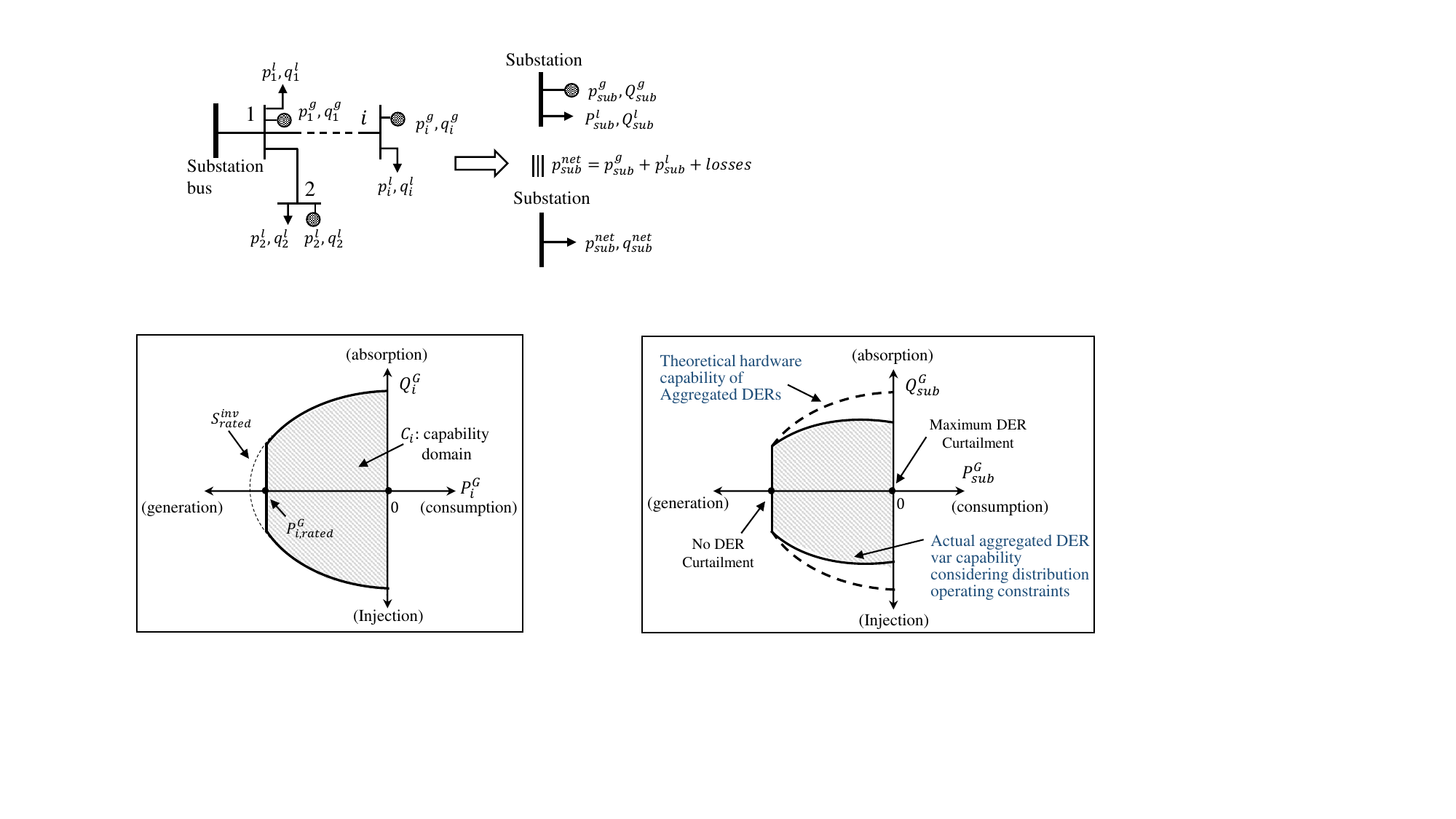}
\caption {One line diagram of a typical distribution feeder with DER and its aggregated representation}
\label{fig:one_line_diagram}
\end{figure}

A typical DS connected to a substation bus with solar PV penetration is shown in \figurename \ref{fig:one_line_diagram}. Load and PV generation at $i^{th}$ node are denoted by $p_i^l+jq_i^l$ and $p_i^g+jq_i^g$ respectively, where $p$ and $q$ denote real and reactive power component respectively. The distribution loads and DERs can be aggregated separately as $p^{l}_{sub}$ and $p^{g}_{sub}$ at the substation as shown in the \figurename \ref{fig:one_line_diagram}. Consequently, the whole DS can further be aggregated as the net power demand at substation which includes actual loads, DERs and losses as shown in the same \figurename \ref{fig:one_line_diagram}. In this section, we will systematically build the characterization of aggregated var capability curve.

\begin{figure}[b]
\centering
\vspace{-4mm}
\includegraphics[trim=-0in 0in 0in 0in,width=2.5in]{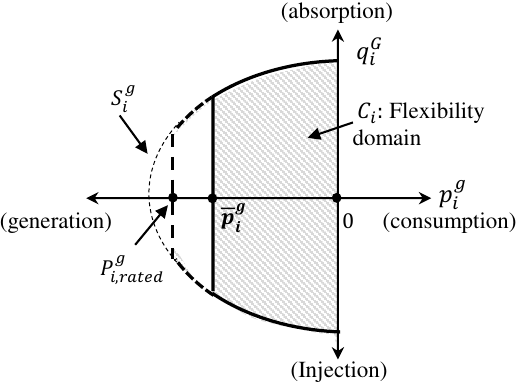}
\caption {Capability curve of a solar PV inverter device }
\label{fig:device_cap_curve}
\end{figure}
\subsection{Var Capability of Individual Solar PV}
For each individual PV inverter, the \textit{device flexibility domain} $\mathcal{C}_i$ can be characterize as following:
\begin{eqnarray}
\mathcal{C}_i=\left \{ (p_i^g,q_i^g) \quad \middle| \quad \begin{array}{l}
	p_i^{g^2}+q_i^{g^2} \leq S_i^{g^2}\\
	p_i^g\leq 0 \\
	|p_i^g|\leq \overline{p}_i^g \leq p^g_{i,rated}
\end{array}
\quad \right \}
\end{eqnarray}
where, $S_i^g$ and $p_{i,rated}^g$ are the hardware capacity of the inverter and solar panel respectively, whereas $\overline{p}_i^g$ is the maximum solar generation possible at given point of time in a day. $\mathcal{C}_i$ represents the available flexibility in var generation or absorption by the inverter for all possible amounts of real power generation. We consider the following sign convention: positive value represents the consumption/absorption and negative value represents generation/injection of real/reactive powers. 
A typical  flexibility domain of a solar PV inverter can be graphically drawn as shown in \figurename \ref{fig:device_cap_curve}.
The outer envelop of the domain $\mathcal{C}_i$ can be defined as a function $q^{g,cap}_i=f(p^g_i)$ which is usually termed as \textit{device Q-capability curve}. This curve is a collection of maximum reactive power values that an individual DER inverter can inject or absorb for a given real power generation. The domain $\mathcal{C}_i$ shrinks or increases as the operating point $\overline{p}^g_i$ moves along the horizontal axis throughout the day.

\subsection{Net DER Aggregation}
\begin{figure}
\centering
\vspace{-4mm}
\includegraphics[trim=-0in 0in 0in 0in,width=2.7in]{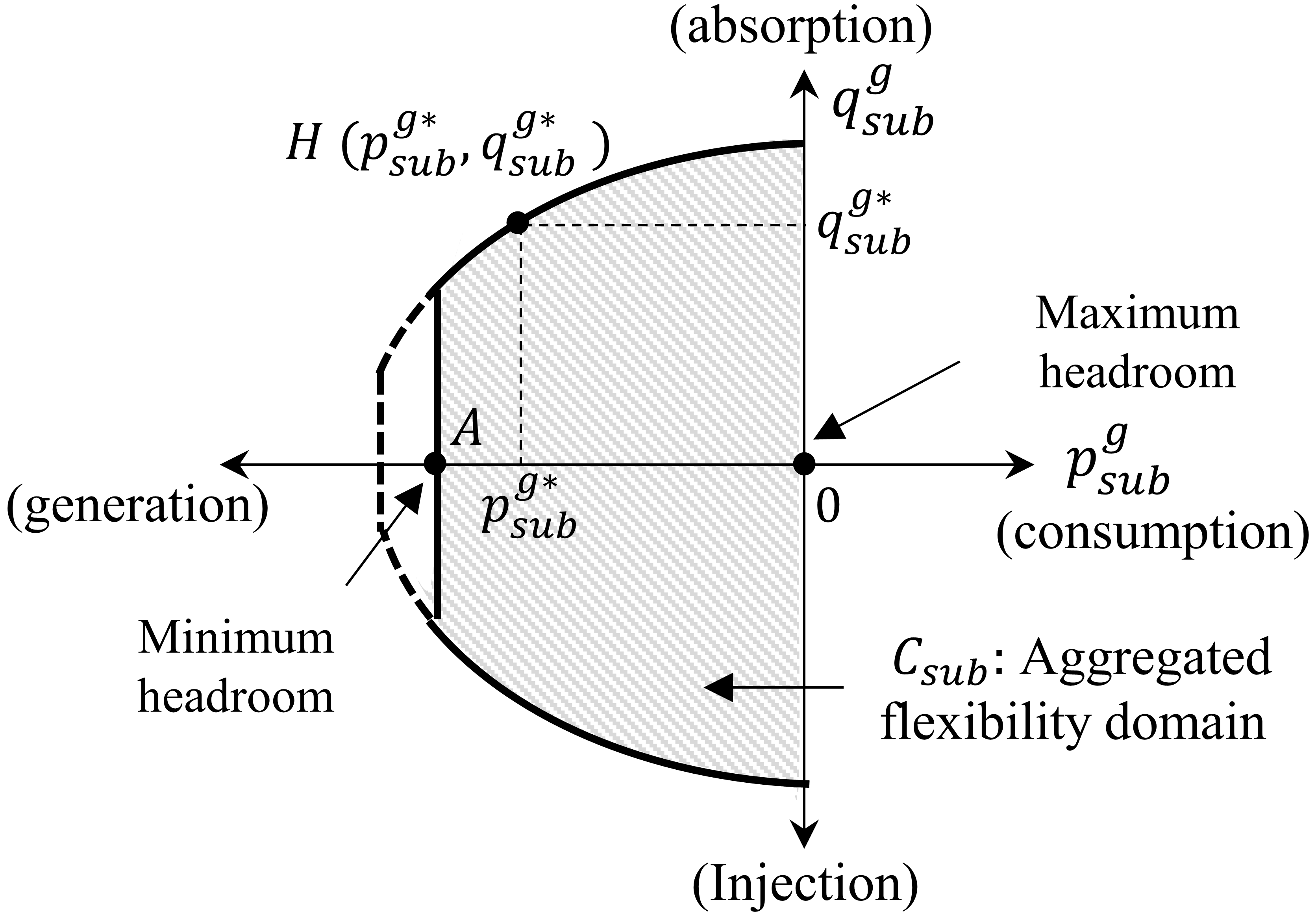}
\vspace{-4mm}
\caption {Aggregated DER Capability curve at Substation}
\label{fig:agg_cap_curve}
\end{figure}

Before developing the net capability of the whole network, lets understand the aggregation of DERs. An \textit{aggregated DER flexibility domain}, $\mathcal{C}_{sub}$, can be defined as the total flexibility provided by all the DERs combined at the substation as following:
\begin{eqnarray}
\mathcal{C}_{sub}=\left \{ (p^g_{sub},q^g_{sub}) \quad \middle| \quad \begin{array}{l}
	p^g_{sub} = \sum_{i=1}^N p_i^g \\
	q^g_{sub} = \sum_{i=1}^N q_i^g \\
	(p_i^{g},q_i^{g}) \in \mathcal{C}_i
\end{array}
\quad \right \}
\end{eqnarray}
where $p^g_{sub}$ and $q^g_{sub}$ are the total real power and var generation from DERs. The outer envelop of the domain $\mathcal{C}_{sub}$ can be defined as a function $q^{cap}_{sub}=f(p^g_{sub})$ that we call as \textit{aggregated DER capability curve} as shown in \figurename \ref{fig:agg_cap_curve}. The horizontal axis can also be seen as varying aggregated real power headroom where point $A$ and origin denote minimum and maximum possible aggregated headroom, respectively. A point $H (p^{g*}_{sub},q^{g*}_{sub})$ on the curve implies that for a given value of $p^{g*}_{sub}$, the maximum possible var absorption is $q^{g*}_{sub}$. Note that a given $p^g_{sub}$ can be achieved in more than one way by different headroom combinations of individual PV generations. In other words, $H$ also denotes the operating point to achieve $q^{g*}_{sub}$ var absorption with minimum possible aggregated headroom. All other possibilities of achieving $q^{g*}_{sub}$ which fall inside the domain will require higher aggregated headroom than necessary. 

\begin{figure}
\centering
\vspace{-4mm}
\includegraphics[trim=-0in 0in 0in 0in,width=3.2in]{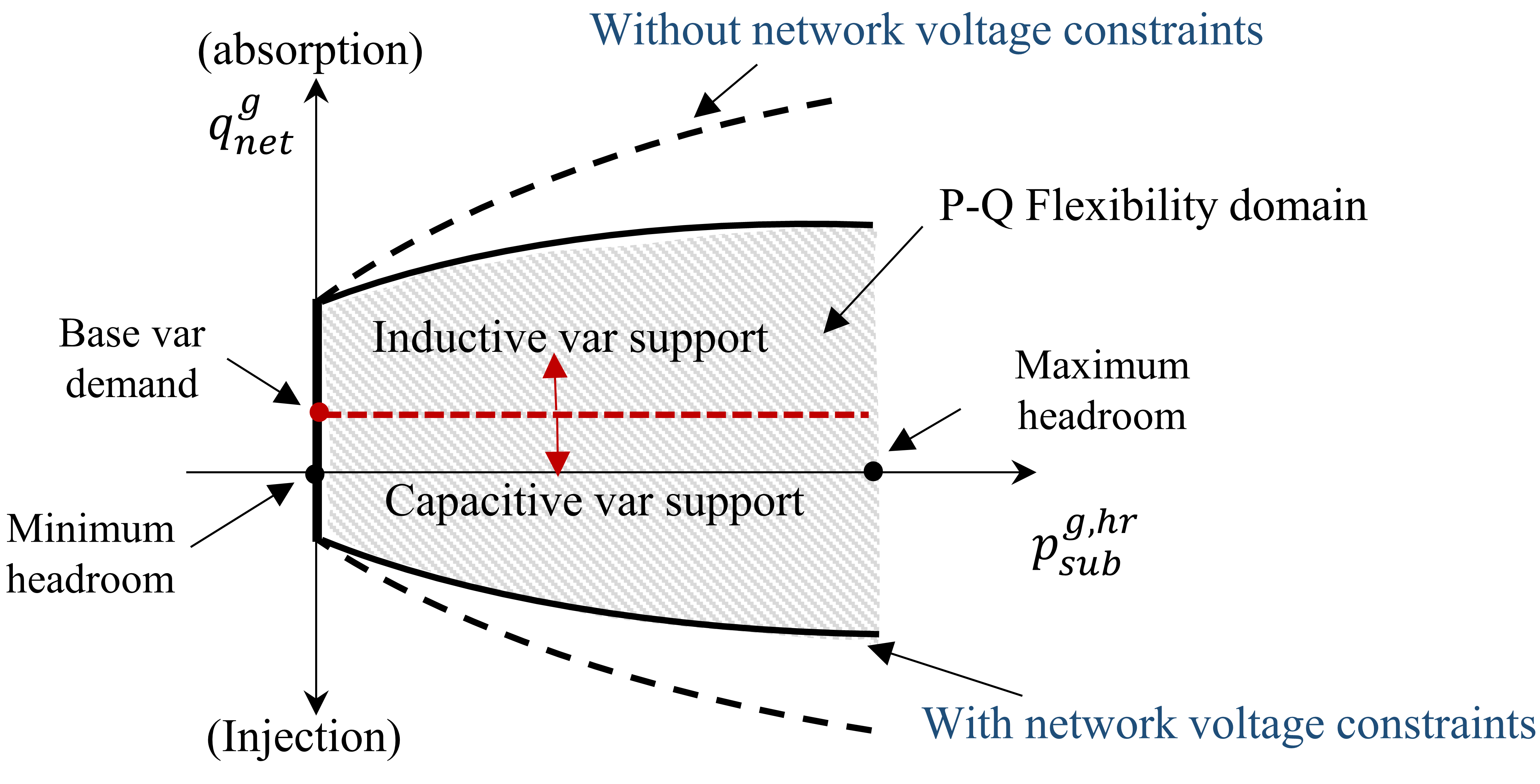}
\caption {Feeder Net Capability curve at Substation}
\label{fig:net_cap_curve}
\end{figure}

However, the more useful information for TSO is the net available var at the substation which includes aggregated load, DER as well as network losses.
Therefore, we define the \textit{aggregated net var capability curve} that provides the information of maximum net var injection/absorption possible at the substation which is seen by the transmission system as net var demand. Henceforth, we will refer it as aggregated capability curve for brevity. We have seen in \figurename \ref{fig:agg_cap_curve} that DER headroom provides real power flexibility that can further enhance the var flexibility region. Therefore, we define aggregated capability curve as function of aggregated DER headroom at the substation, $p^{g,hr}_{sub}$, i.e. $q_{net}^{cap}=f(p^{g,hr}_{sub})$. A conceptual curve at a given operating condition is shown in \figurename \ref{fig:net_cap_curve} that depicts the capacitive and inductive var flexibility domain.
Another important point to consider is that the inverter var injection or absorption affects the voltage profile of the DS and consideration of voltage limits may shrink the flexibility domain in certain operating conditions  as visible in the \figurename \ref{fig:net_cap_curve}. 

The total headroom, $p^{hr}_{sub}$, aggregated at the substation level can be defined such that
{\color{blue}
\begin{equation}
p^{g}_{sub} = \sum_i{S}^{g}_{i} .(1-p^{g,hr}_{sub}) 
\end{equation}
$p^{g,hr}_{sub}$ ranges from $p^{hr}_{min}$ to 1, where, $p^{hr}_{min}$ is the minimum possible aggregated headroom at a give time in a day. $p^{hr}_{min}$ corresponds to the maximum generation possible at given time.
\begin{equation}
\sum_i{S}^{g}_{i} .(1-p^{hr}_{min}) = \sum_i \overline{p}^{g}_{i} 
\end{equation}
Thus, $p^{hr}_{min}$ can be obtained as,
\begin{equation}
p^{hr}_{min} = 1 - \frac{\sum_i \overline{p}^{g}_{i}}{\sum_i{S}^{g}_{i}} 
\end{equation}
Now, for each DER, we can write,
\begin{equation}
p^{g}_{i} = S^{g}_{i}.(1-p^{g,hr}_{i})
\end{equation}
where $p^{g,hr}_{i}$ is headroom for $i^{th}$ DER. Finally, $p^{g,hr}_{sub}$ can be written in form of $S_i^g$ and $p_i^{g,hr}$ as
\begin{equation}
p^{g,hr}_{sub} = (\sum_{i=1}^N {S}^{g}_{i}.p^{g,hr}_{i})/\sum_i S^{g}_{i}
\end{equation}
}
The real power headroom for DER can be achieved via inverter oversize, storage or curtailment. It exhibits the higher flexibility of the system and provide more options to TSO to handle var related grid events. Nonetheless, utilizing this flexibility involves a greater discussion on policy, customer comfort, and related cost-benefit analysis which is beyond the scope of this paper.

\section{Process of Capability Estimation}
\subsection{System Modeling}
In this section, we will utilize the linearized \textit{'LinDist3 Flow'} equations \cite{arnold_optimal_2016} for an unbalanced three-phase DS to develop a graph-representation model \cite{zhu_fast_2016}. Consider a radial DS with $N+1$ nodes represented by a tree graph $\mathcal{T}=(\mathcal{N},\mathcal{E}$), where $\mathcal{N}:=\{0,1,\cdots,N\}$ is a set of DS nodes, indexed by $i$ and $j$. For simplicity, let's assume each $i\in \mathcal{N}$ has all three phases $a,b$ and $c$. The set $\mathcal{E}:=\{(i,j)\}$ contains all line segments with $i$ as the upstream and $j$ as the downstream node. Each line element $(i,j) \in \mathcal{E}$ will also have three phases. The subset $\mathcal{N}_j$ is a collection of all immediate downstream neighboring buses of node $j$. The secondary side of the substation is denoted by node 0. 
Let $M$ be an $3N\times 3N$ graph incidence matrix of $\mathcal{T}$. The $l^{th}$ column of matrix $M$ corresponds to the line segment $(i,j)\in \mathcal{E}$ with entries $M(i,k)=e$ and $M(j,k)=-e$, where $e$ is a $3 \times 3$ identity matrix. All other entries of $M$ are zero. Now, according to \textit{LinDist3Flow} model, the voltages at node $i$ and $j$ can be written as
\begin{equation}
\label{eq:LinDist1}
\mathbb{V}_i\mathbb{V}_i^* = \mathbb{V}_j \mathbb{V}_j^* - \mathbb{Z}_{ij}^p\mathbb{P}_j - \mathbb{Z}_{ij}^q\mathbb{Q}_j 
\end{equation}
where, $\mathbb{V}_j=[V_a V_b V_c]_j^T$ represent the vector of voltage phasors at node $j$. Similarly, $\mathbb{P}_j\!=\![P_a \ P_b \ P_c]_j^T$ and 
$\mathbb{Q}_j\!=\! [Q_a \ Q_b \ Q_c]_j^T$ denote the real and reactive power entering at node $j$. $\mathbb{Z}_{ij}^q$ and $\mathbb{Z}_{ij}^p$ are the constant three phase impedance matrices for line segment $(i,j)$ as following \cite{arnold_optimal_2016}:

\begin{equation}
\color{blue}
    \mathbb{Z}_{ij}^P = \begin{bmatrix}
                        -2r_{aa} & r_{ab}-\sqrt{3}x_{ab} & r_{ac}+\sqrt{3}x_{ac} \\
                        r_{ba}+\sqrt{3}x_{ba} & -2r_{bb} & r_{bc}-\sqrt{3}x_{bc} \\
                        r_{ca}-\sqrt{3}x_{ca} & r_{cb}+\sqrt{3}x_{cb} & -2r_{cc}
                    
                    \end{bmatrix}_{ij}
\end{equation}
\begin{equation}
\color{blue}
    \mathbb{Z}_{ij}^Q = \begin{bmatrix}
                        -2x_{aa} & x_{ab}+\sqrt{3}r_{ab} & x_{ac}-\sqrt{3}r_{ac} \\
                        x_{ba}-\sqrt{3}r_{ba} & -2x_{bb} & x_{bc}+\sqrt{3}r_{bc} \\
                        x_{ca}+\sqrt{3}r_{ca} & x_{cb}-\sqrt{3}r_{cb} & -2x_{cc}
                    
                    \end{bmatrix}_{ij}
\end{equation}
where, $Z_{\phi\psi, ij} =  r_{\phi\psi, ij} + j x_{\phi\psi, ij}$ is impedance between phase $\phi$ and $\psi$ of line segment (i,j).
Now, let's define the vector of squared of voltage magnitude as a new variable $\mathbb{Y}_j\! =\! \mathbb{V}_j\mathbb{V}_j^*\!=\![y_a \ y_b \ y_c]_j$ for $j \in \mathcal{N} \backslash \{0\}$. Assuming the reference node 0 voltage as $\mathbb{Y}_0$, the voltages at each node can be written in compact form as following:
\begin{equation}
\label{eq:LinDist2}
[M_0 \ M^T][\mathbb{Y}_0^T \ \mathbb{Y}^T]^T= M_0\mathbb{Y}_0+M^T\mathbb{Y}=- \mathbb{Z}_{D}^p\mathbb{P} - \mathbb{Z}_{D}^q\mathbb{Q} 
\end{equation}
Where, $M_0$ is a matrix of size $3N \times 3$ with first entry as $e$ and rest as zero. $\mathbb{Z}_{D}^p$ and $\mathbb{Z}_{D}^q$ are diagonal matrices of size $N$ where $l^{th}$ entries are $\mathbb{Z}_{ij}^p$ and $\mathbb{Z}_{ij}^q$ respectively which correspond to $l^{th}$ line segment $(i,j)$. 

The line flows $\mathbb{S}_j =\mathbb{P}_j + j\mathbb{Q}_j$ can be written in form of net injections as following:


\begin{equation}
\label{eq:LinDist3}
\mathbb{S}_j\approx -s_j+\sum_{k \in \mathcal{N}_j}\mathbb{S}_k + L_j
\end{equation}
Where $s_j=p_j+q_j$ is the vector of net injection at node $j$ at all phases denoted by $s_{\phi,j}$ where, $\phi \in {a,b,c}$. Usually, in LinDistFlow model, line losses are neglected which introduce a relatively small error in the modeling as indicated by \cite{farivar_equilibrium_2013}. However, to increase accuracy, we consider a constant loss term $L_j$ in (\ref{eq:LinDist3}). The loss term $L_j$ denotes the losses incurred in line ending at node $j$ and can be estimated based on the offline study of the base operating point as indicated in \cite{zhu_fast_2016}. Equation (\ref{eq:LinDist3}) can be re-written in compact form as
\begin{align}
\label{eq:LinDist4a}
-M\mathbb{P} = -p + L_p \\
\label{eq:LinDist4b}
-M\mathbb{Q} = -q + L_q
\end{align}
Where $L_p$ and $L_q$ are vectors of real and reactive loss factors. Using (\ref{eq:LinDist2}), (\ref{eq:LinDist4a}) and (\ref{eq:LinDist4b}), voltages in form of net injections can be written as following:

\begin{equation}
\label{eq:LinDist5}
\mathbb{Y}=R^{eq}p +X^{eq}q - M^{-T}M_0\mathbb{Y}_0 + L_c
\end{equation}
Where $R^{eq}=-M^{-T}\mathbb{Z}_{D}^p M^{-1}$, $X^{eq}=-M^{-T}\mathbb{Z}_{D}^q M^{-1}$ and $L_c = R^{eq}L_p-X^{eq}L_q$ are constants. Lets assume that the substation voltage is balanced and has same magnitude in all phases denoted by a a scalar $v_0$. Further, due to radial structure of network, $M^{-T}M_0\mathbb{Y}_0$ is same as -$v_0^2\mathbb{I}$, where $\mathbb{I}$ is a column vector of size $3N$ with all entries as 1.  

Substation secondary voltage $v_0$ can be controlled via an on-load tap changer (OLTC) within a range as $v_0 = {v^{tm}.r} $, where $v^{tm}$ is primary side transmission voltage and $r$ is tap ratio of OLTC. Usually each tap provides $\pm 0.0063$ pu voltage regulation with maximum $\pm 16$ taps. Therefore the maximum possible values of $r$ are $1\pm 0.1$.


Let's assume the DERs are located at the nodes collected in a subset $\mathcal{G}\subseteq \mathcal{N}$. In this case, only inverter-based DERs are considered such as solar PV. The net power injection of real and reactive power at each node $j \in \mathcal{G}$ is denoted by $p_{\phi,j}=-p_{\phi,j}^g - p^l_{\phi,j}$ and $q_{\phi,j}=-q_{\phi,j}^{g} - q^l_{\phi,j}$ respectively. Superscript $g$ and $l$ denote DER and loads. 
All loads are constant power loads and the capacitors are modeled as reactive power loads. 

\color{blue}
\subsection{Estimating System Losses and Net Var Demand}
In a multi-phase unbalanced system, losses also include the impact of mutual coupling impedance with other phases. In a 3-phase unbalanced system, losses incurred in each phase of the line ending at node $j$, $\mathcal{L}_j$, can be written in vector form as

\begin{equation}
\label{eq:losses}
\mathcal{L}_j = d(\mathbb{I}_j\mathbb{Z}_j\mathbb{I}_j^T) = \begin{bmatrix}
         Z_{aa}I_{a}I_{a}^*+Z_{ab}I_{b}I_{a}^*+Z_{ac}I_{c}I_{a}^* \\
         Z_{ba}I_{a}I_{b}^*+Z_{bb}I_{b}I_{b}^*+Z_{bc}I_{c}I_{b}^* \\
          Z_{ca}I_{a}I_{c}^*+Z_{cb}I_{b}I_{c}^*+Z_{cc}I_{c}I_{c}^*
         \end{bmatrix}_{j}
\end{equation}
where, $\mathbb{I}_j = [I_{a} I_{b} I_{c}]_{j}$ and $\mathbb{Z}_j$ are current and impedance of line ending at node $j$. $Z_aa$ and $Z_{ab}$ represent self impedance of phase a and mutual impedance between phase a and b, and likewise. $d(X)$ returns the diagonals of a square matrix X as a vector.  Following \cite{schweitzer_lossy_2020,low_convex_relaxation}, voltage phasors are assumed to be balanced, i.e.,  $V_{a,j}V_{b,j}^{-1} \approx V_{b,j}V_{c,j}^{-1} \approx V_{c,j}V_{a,j}^{-1} \approx \alpha = e^{j2\pi/3}$. Note that $1/\alpha = \alpha^2$.
Using this approximation, currents in each phase can be written as 
\begin{align}
\label{eq:Ia_Sa}
\begin{split}
I_a^* = S_a/V_a = \alpha^2 S_a/V_b = \alpha S_a/V_c\\
I_b^* = \alpha S_b/V_a = S_b/V_b  = \alpha^2 S_c/V_c\\
I_c^* = \alpha^2 S_c/V_a = \alpha S_c/V_b =S_c/V_c
\end{split}
\end{align}
Replacing (\ref{eq:Ia_Sa}) in (\ref{eq:losses}) such that losses in phase $\phi$ is function of $V_{\phi}$, the loss term can be written as
\begin{equation}
\label{eq:loss_term}
     \mathcal{L}_j= \begin{bmatrix}
         S_{a}/y_{a}\\
         S_{b}/y_{b}\\
        S_{c}/y_{c}
         \end{bmatrix}_{j}   
         \begin{bmatrix}
           Z_{aa} & Z_{ab}\alpha &
         Z_{ac}\alpha^2\\
        Z_{ba}\alpha^2 & Z_{bb} &
         Z_{bc}\alpha\\
         Z_{ca}\alpha & Z_{cb}\alpha^2 &
         Z_{cc}\\
         \end{bmatrix}_{j}
         \begin{bmatrix}
         S_{a}^*\\
         S_{b}^*\\
         S_{c}^*
         \end{bmatrix}_{j} 
\end{equation}
Where, ${y_{\phi,j}}={V_{\phi,j}^2}$. Henceforth, we drop the subscript $\phi$ for the convenience of the notations. Equation (\ref{eq:loss_term}) is able to capture the impact of mutual impedance of other phases as well which may be significant in an unbalanced distribution systems.

The net reactive power demand at the substation, $q_{sub}^{net}$ 
can be written as sum of the net injection of var at each node due to load, DER inverter and reactive power losses incurred at each line across all three phases as following.
\begin{equation}
\small
\label{eq:q_net}
q_{sub}^{net}(q^g_j,y_j^g) =\sum_{\phi \in \{a,b,c\}} \Big(\sum_{j\in \mathcal{G}}q_{\phi,j}^{g}-\sum_{j\in \mathcal{N}}q_{\phi,j}^{l} + \sum_{j\in \mathcal{N}}\mathcal{L}_{\phi,j}\Big)
\end{equation}

\color{black}
\subsection{DER-OPF Formulation}
Our objective here is to construct the net capability curve $q_{net}^{cap}=f(p^{g,hr}_{sub})$ as shown in \figurename \ref{fig:net_cap_curve}. To achieve it, we need to estimate both the capacitive ($\underline{q}_{net}^{cap}$) and inductive ($\overline{q}^{cap}_{net}$) var capabilities of the network which is same as minimizing and maximizing the net var flow at the substation. Based on the already defined preliminaries, following DER-OPF can be written to estimate $\underline{q}_{net}^{cap}$ :

\begin{mini!}|l|[3]
{\small p_j^{g,hr},q_{j}^{g},v_0}{q^{net}_{sub}(y_{j},q_{j}^{g})}
{\label{opt}}{}
\vspace{-1mm}
\addConstraint{\mathbb{Y}=R^{eq}p +X^{eq}q +({v_0})^2\mathbb{I} + L_c}{\label{opt:lindist}}{\forall j\in \mathcal{N}}
\addConstraint{p_{j}=S_{j}^g.(1-p_j^{g,hr}) - p^l_{j}}{,\label{opt:DGP_const}\quad}{\forall j\in \mathcal{N}}
\addConstraint{q_{j}=q_{j}^{g} - q^l_{j}}{,\label{opt:DGQ_const}\quad}{\forall j\in \mathcal{N}}
\addConstraint{\underline{y}\leq y_{j}\leq \overline{y}}{,\label{opt:volt_const}\quad}{\forall j\in \mathcal{N}}
\addConstraint{|q_{j}^{g}|\leq \sqrt{S^{g^2}_j\!\!-S_{j}^{g^2}(1\!-p_j^{g,hr})^2} }{\label{opt:q_const}\quad}{\forall j\in \mathcal{G}}
\addConstraint{p_{min}^{hr}\leq p_j^{g,hr}\leq 1}{,\label{opt:curt_const}\quad}{\forall j\in \mathcal{G}}
\addConstraint{\sum_{j\in \mathcal{G}} S^{g}_{j}.p^{g,hr}_{j} = p^{g,hr}_{sub}.\sum_{j\in \mathcal{G}}S^{g}_{j} }{\label{opt:curt_sum}\quad}{}
\addConstraint{v_{tm}\underline{r}\leq v_0 \leq v_{tm}\overline{r}}{,\label{opt:tap_const}\quad}{}
\end{mini!}
 Constraints (\ref{opt:lindist} - \ref{opt:DGQ_const}) denote the power flow and (\ref{opt:volt_const}) ensures the voltages are within the ANSI limits \cite{noauthor_ansi_2016}. $\overline{y}$ and $\underline{y}$ are upper and lower allowable voltage limits, and are usually taken as {\color{blue} $1.05$ and $0.95$, respectively}. Constraint (\ref{opt:q_const}) manifest the hardware capacity limit of an inverter. Constraint 
(\ref{opt:curt_const}) provides the range of headroom. Note that the aggregated headroom at the substation level, $p^{g,hr}_{sub}$, is an input parameter in this optimization whereas individual DER headroom, $p_j^{g,hr}$, is an optimization variable. The relation between these two is given by (\ref{opt:curt_sum}). To avoid integer programming, $r$ is taken as a continuous variable in the OLTC constraint in (\ref{opt:tap_const}). Upper and lower saturation limits on OLTC tap ratios are denoted by $\overline{r}$ and $\underline{r}$ respectively. \color{blue}Note that the $v_{tm}$ is an input parameter to the formulation whereas $v_0$ is an optimization variable. This allow us to perform parametric study of capability curve with respect to grid side voltage.
\color{black}
The solution of the optimization (\ref{opt}) provides the optimal var set dispatch ($q_{j}^{g*}$) and optimal headroom ($p_{j}^{g,hr*}$) for each DER and optimal secondary side voltage set-point ($v_0^*$).

Similar to (\ref{opt}), upper part of the net capability curve $(\overline{q}_{net}^{cap})$ can be estimated by maximizing the net var demand at substation which is same  as following,
\begin{mini!}|l|[3]
{\small p_j^{g,hr},q_{j}^{g},y_{0}}{-q^{net}_{sub}(y_{j},q_{j}^{g})}
{\label{opt1}}{}
\vspace{-1mm}
\addConstraint{(\ref{opt:lindist})-(\ref{opt:tap_const})}{}{\forall j\in \mathcal{N}}
\end{mini!}
Unfortunately, the objective function in (\ref{opt1}) is not convex due to losses term $\mathcal{L}_j$ in (\ref{eq:q_net}) being quadratic as shown in (\ref{eq:loss_term}). However, it can be converted to a convex expression by removing the $\mathcal{L}_j$ term as following:
\begin{mini!}|l|[3]
{\small p_j^{g,hr},q_{j}^{g},y_{0}}{- \Big(\sum_{j\in \mathcal{G}}q_{\phi,j}^{g}-\sum_{j\in \mathcal{N}}q_{\phi,j}^{l}\Big)}
{\label{opt2}}{}
\vspace{-1mm}
\addConstraint{(\ref{opt:lindist})-(\ref{opt:tap_const})}{}{\forall j\in \mathcal{N}}
\end{mini!}
Usually, the var losses are much smaller component of $q_{sub}^{net}$ compared to combined var consumption by the loads and the inverters, therefore, it doesn't affect the optimal point significantly. In fact, in most cases, the optimal point of (\ref{opt2}) is also optimal for (\ref{opt1}) except when lower voltage boundary constraints of (\ref{opt:volt_const}) at all nodes are not active. In those cases, (\ref{opt1}) tries to further reduce voltage to its minimum in order to increase losses which adds a negligible error in optimal net var flow $q_{sub}^{net}$ calculated by (\ref{opt2}). Therefore, $\overline{q}_{cap}^{net}$ is estimated via (\ref{eq:q_net}) using optimal $q^{g*}_j$ resulting from (\ref{opt2}). 

\color{blue}
\subsection{Decoupling Range of Capability Curve with TN Voltage}
It is pertinent to discuss that 
OLTC tap ratio provides a limited decoupling between primary and secondary side of the substation within the range of $r$. Due to this, the desired optimal secondary voltage $v^*_0$ can be achieved by adjusting tap ratio for any value of $v_{tm}$ which lies in the decoupling range $\mathcal{D}$ defined as,
\begin{equation}
\mathcal{D}=[\:v^*_0/\overline{r}, \; v^*_0/\underline{r}\:]
\end{equation}
This decoupling is lost when $v_{tm} \notin \mathcal{D}$ i.e. the OLTC tap gets saturated. 
To address this concern, DSO estimates the var capability $q_{net}^{cap}=[ \underline{q}_{net}^{cap}, \;\; \overline{q}_{net}^{cap} ]$ and $\mathcal{D}$ for nominal value of $v_{tm}=1$ and send this information to TSO to be included in their optimizations. 
It is advised that the TSO use both $q_{net}^{cap}$ and $\mathcal{D}$ as constraints while estimating the var requirement service from DSO. However, there might arise situations when expected $v_{tm} \notin \mathcal{D}$. In such cases, DSO can not guarantee the $q_{net}^{cap}$ but can provide an estimated bound on $q_{net}^{cap}$ for worst case value of $v_{tm}$ i.e. $q_{net}^{cap}=[ \underline{q}_{net}^{cap}-\underline{\epsilon}_q, \;\; \overline{q}_{net}^{cap}-\overline{\epsilon}_q ]$ by performing paramteric study of the capability curve estimation with respect to $v_{tm}$ as also shown in the results section later. 
\begin{figure}
\centering
\includegraphics[trim=-0in 0in 0in 0in,width=3.5in]{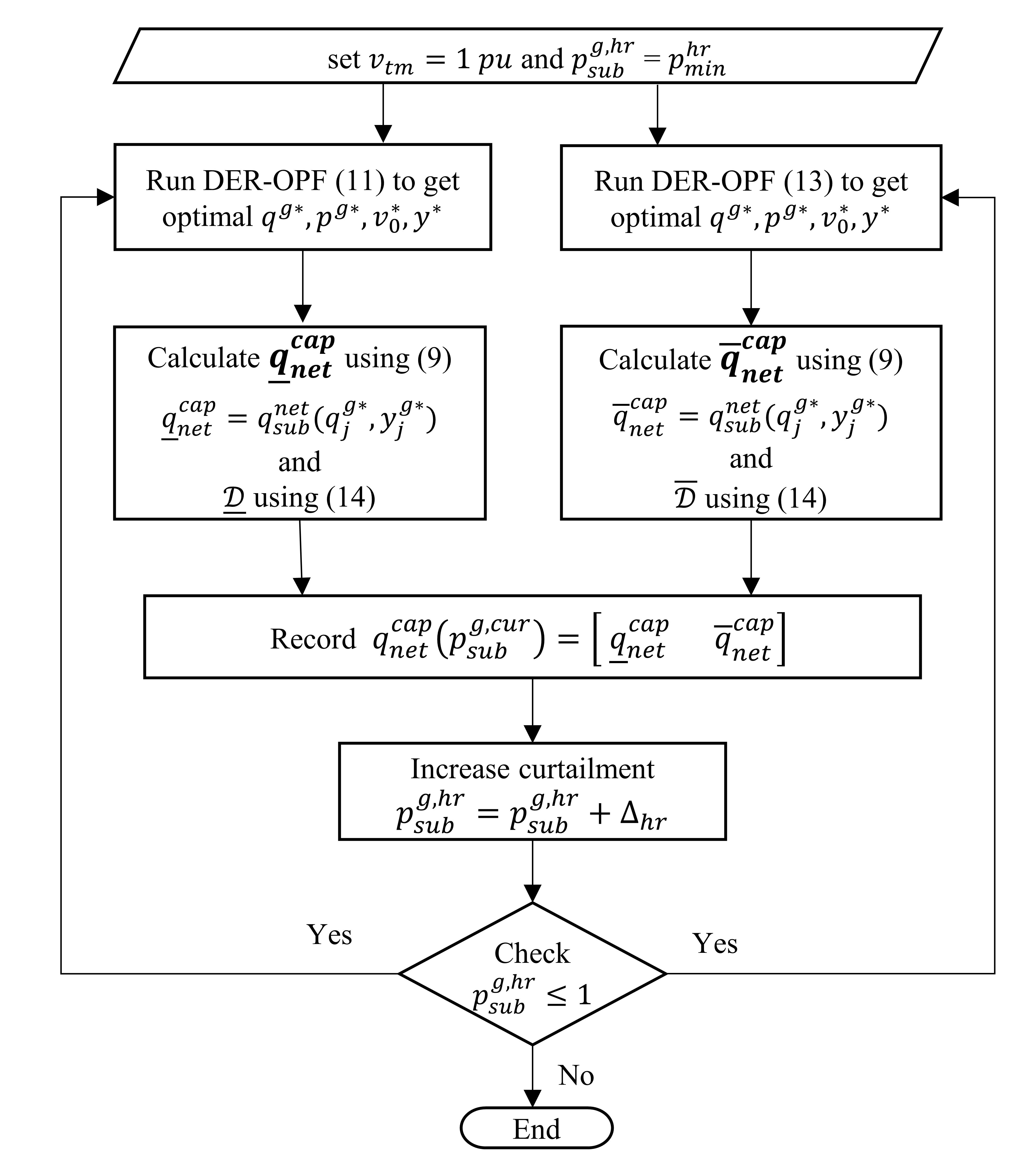}
\vspace{-3mm}
\caption {Flow chart of the process of estimating var capability curve as function of DER Headroom}
\label{fig:FlowChart}
\end{figure}
\color{black}
Flow chart of the overall process of the var capability curve estimation for a given operating condition is shown in \figurename \ref{fig:FlowChart}. 

\vspace{-1mm}


\section{Test Case Study}
\subsection{Reactive Power Flexibility Region (RPFR)}
In order to numerically evaluate the var flexibility provided by DER, we define reactive power flexibility region (RPFR) as the range $[a, \; b]$ at any give operating condition, 
\begin{equation}
a=\underline{q}_{net}^{cap}-q^{base}_{net}\;\;\;\;\;\;\;\; b=\overline{q}_{net}^{cap}-q^{base}_{net}
\label{eq:RPFR}
\end{equation}
where $q^{base}_{net}$ is the net var demand at the substation when all DER inverters are operating in unity power factor mode. 
The $a$ and $b$ denote the maximum available capacitive and inductive var support in MVar, respectively. A higher magnitudes of both $a$ and $b$ with negative and positive signs respectively represent a larger flexibility region. A zero value of both $a$ and $b$ denotes no available var flexibility.  

\subsection{ Test System Description}
\begin{figure}
\centering
\subfloat[][IEEE 37 node test system]{\includegraphics[clip,width=0.7\columnwidth]{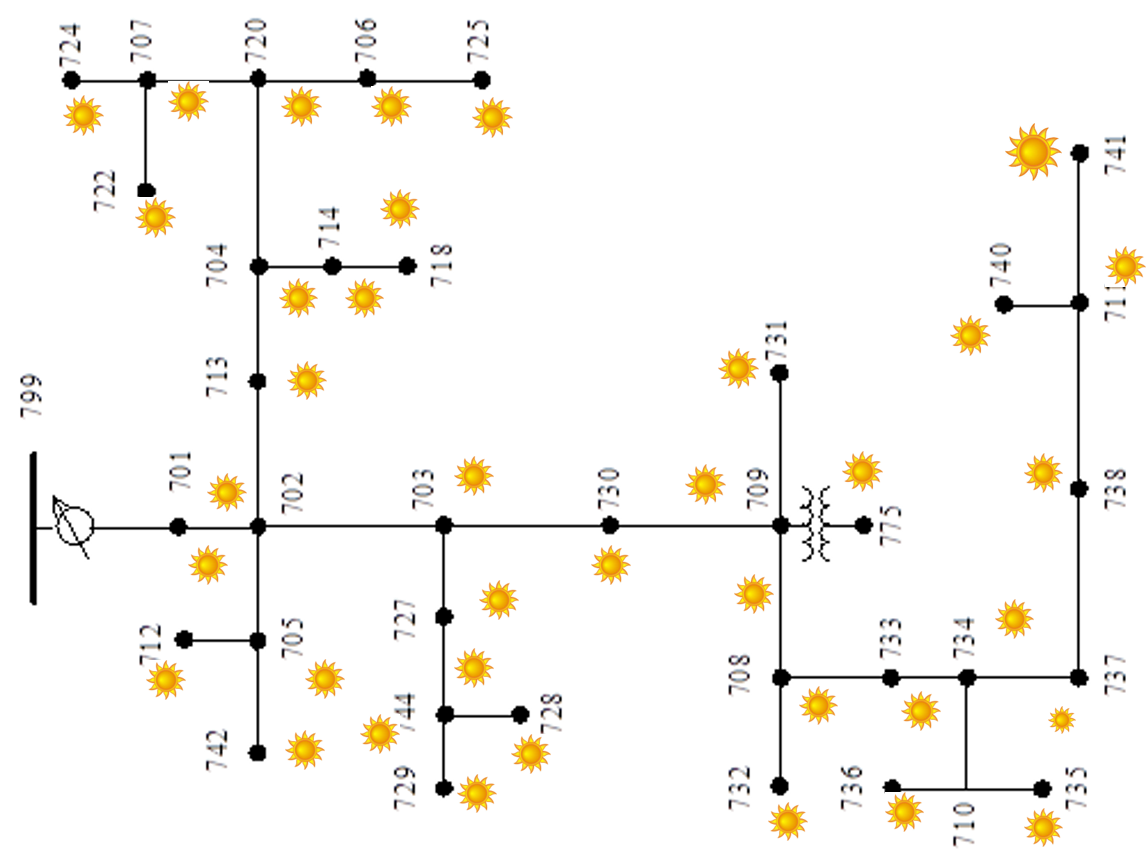}\label{fig:37DSystem}}
\newline
\subfloat[][IEEE 123 node test system]{\includegraphics[clip,width=0.8\columnwidth]{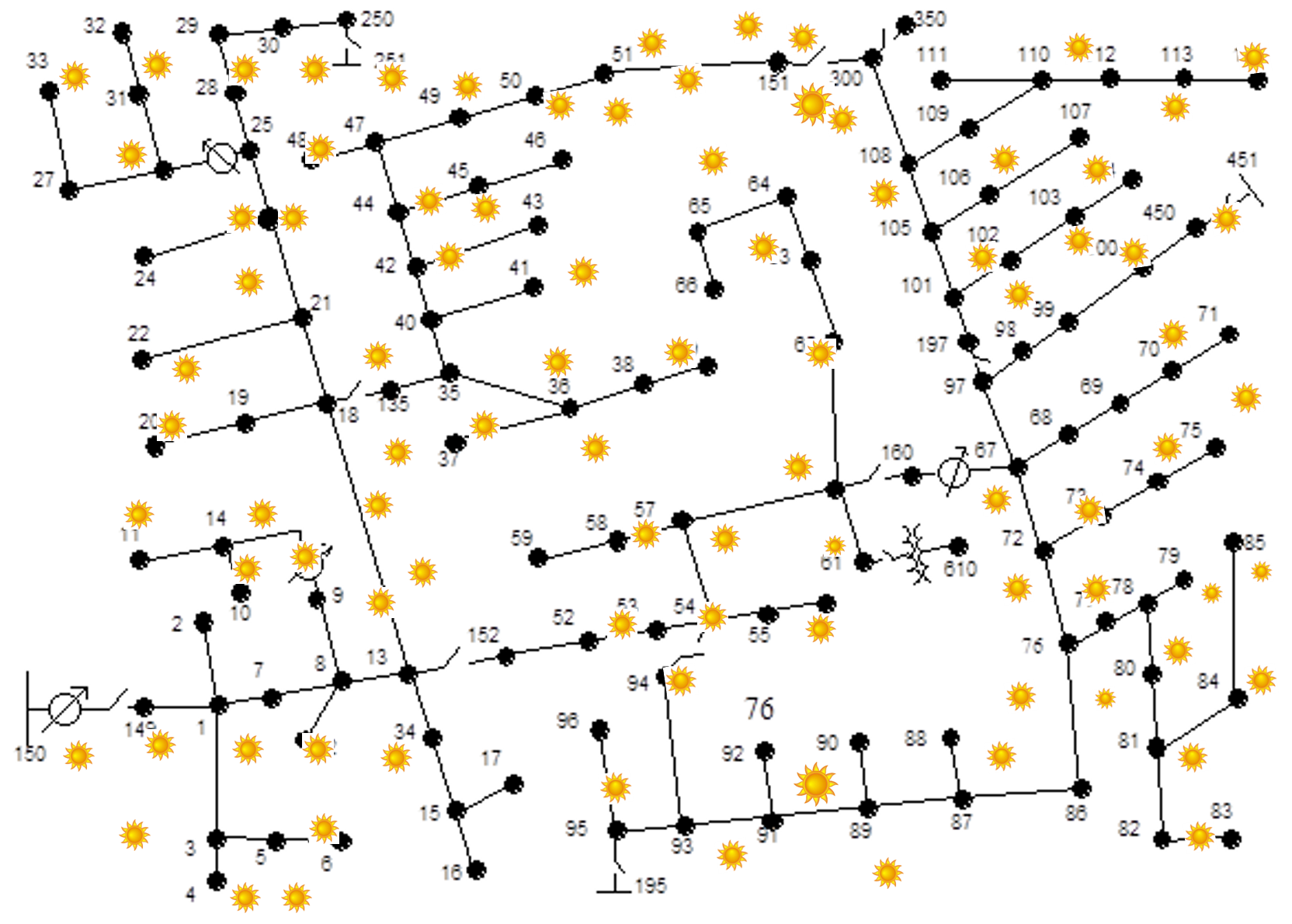}\label{fig:123DSystem}}
\vspace{-2mm}
\caption {Unbalanced distribution test feeders}
\vspace{-1mm}
\label{fig:TestDSystem}
\end{figure}

An unbalanced 3-phase IEEE distribution 37 bus test system is considered with around 4 MW as peak load and around 90\% solar PV penetration as shown in \figurename \ref{fig:37DSystem}. Here, we define the penetration level is a ratio of peak solar generation to peak load demand. Around 100 Single phase DER (solar PV) units are equally distributed throughout the DS nodes in all three phases. Inverter ratings are considered 1.1 times the peak solar generation. Maximum and minimum values of $v_{tm}$ are considered as 0.9 and 1.1.

\color{blue}
\subsection{Modeling Validation}
\begin{figure}
    \color{blue}
	\centering
	\includegraphics[trim=0in 0in 0in 0in,width=1\columnwidth]{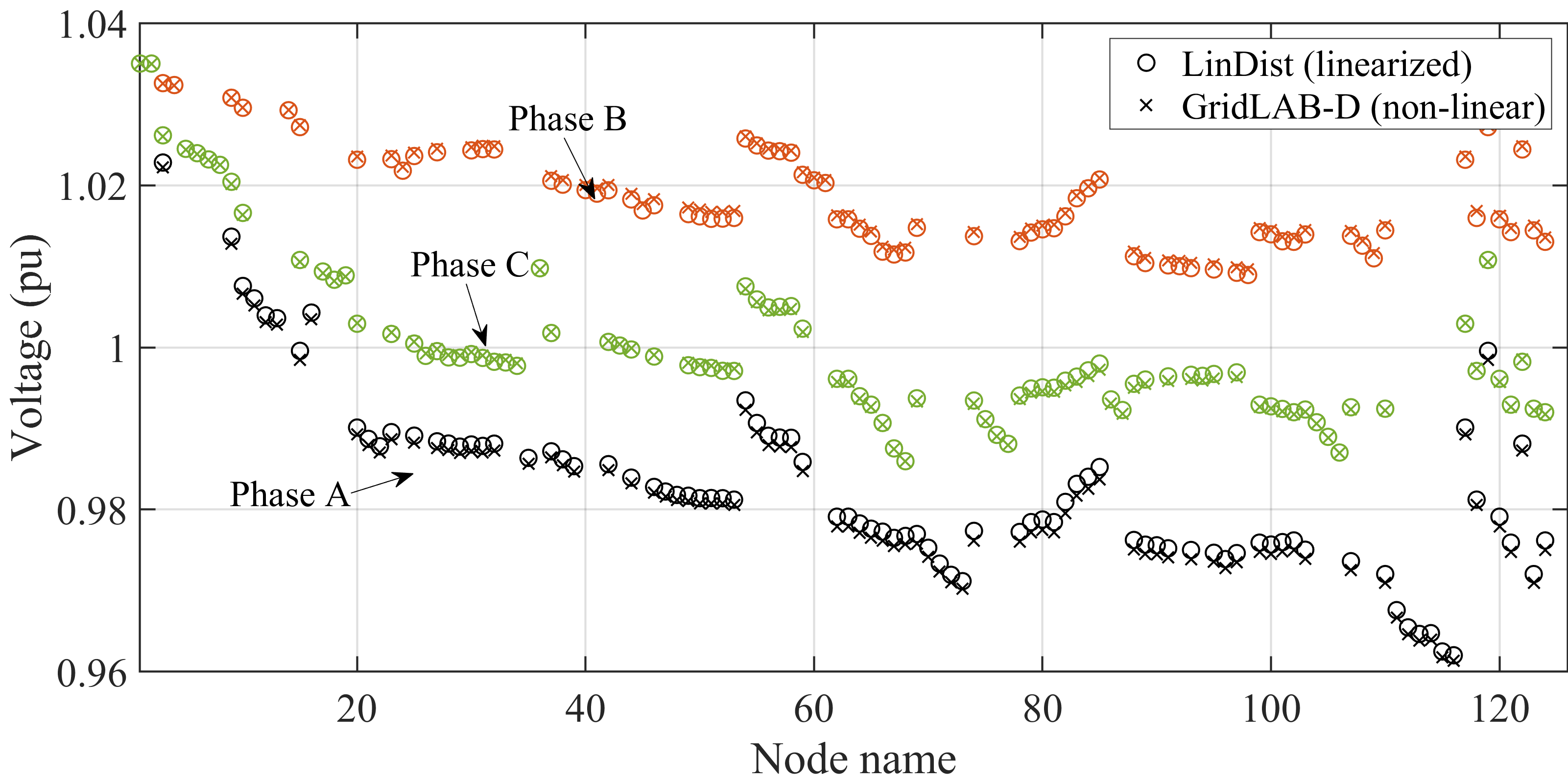}
    \caption {Comparison of voltage profiles obtained from linearized model used in OPF and full non-linear model (GridlabD) to check accuracy}
    \label{fig:volt_compare}
\end{figure}

The successful implementation of the proposed var support framework is highly dependent on the accuracy of the \textit{LinDist3Flow} model used for unbalanced distribution feeder. To verify, the voltages profiles from this model are compared with the full nonlinear model of the DS (GridlabD) for both IEEE 37 and 123 test systems. \figurename \ref{fig:volt_compare} compares the voltages of all three phases of 123 bus system at all the nodes at peak load scenario without any DERs. It can be seen that the voltage profile from linearized model is very close to the non-linear solver. Table \ref{tab:lin_dist_validation} shows that the mean and maximum voltage errors are within 0.02\% and 0.07\% respectively. The high accuracy of LinDist3 model in base case is due to an appropriate selection of loss factor based on offline study. We will keep the same loss factor for capability curve estimations and report the validation error in next subsections. Note that the error in substation net var for both systems is within 0.02\% and 1\% respectively. 
\begin{table}[]
\color{blue}
\caption{Accuracy of LinDist model compared to the non-linear model in calculating voltages and substation reactive power in base case with no DER  }
\label{tab:lin_dist_validation}
\centering
\renewcommand{\arraystretch}{1.3}
\begin{tabular}{cccc}
\hline
\multirow{2}{*}{Test System} & \multicolumn{2}{l}{Voltage error (\%)} & 
\multicolumn{1}{c}{\multirow{2}{*}{\begin{tabular}[c]{@{}c@{}}Substation var \\ error (\%)\end{tabular}}} \\ \cline{2-3}
 & \multicolumn{1}{c}{Mean} & \multicolumn{1}{c}{Maximum} & \multicolumn{1}{c}{}\\
 \hline
IEEE 37 & 0.01 & 0.05 & 0.16 \\
IEEE 123 & 0.02 & 0.07 & 0.73  \\ \hline
\end{tabular}
\end{table}

\color{black}
\subsection{Aggregated Net Capability Curves}
\begin{figure}
\centering
\includegraphics[trim=-0in 0in 0in 0in,width=3.2in]{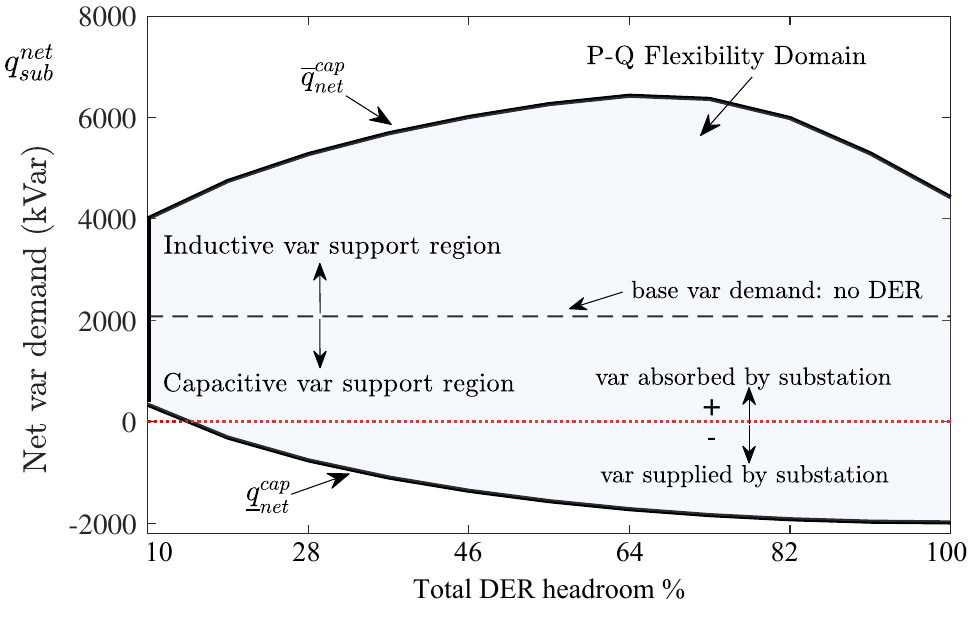}
\vspace{-2mm}
\caption {The aggregated net capability curve of IEEE 37 node DS with high DER penetration as function of DER headroom for case 1 (peak load)} 
\vspace{-0mm}
\label{fig:cap_highload}
\end{figure}
\begin{figure}
\centering
\vspace{-4mm}
\includegraphics[trim=-0in 0in 0in 0in,width=3.2in]{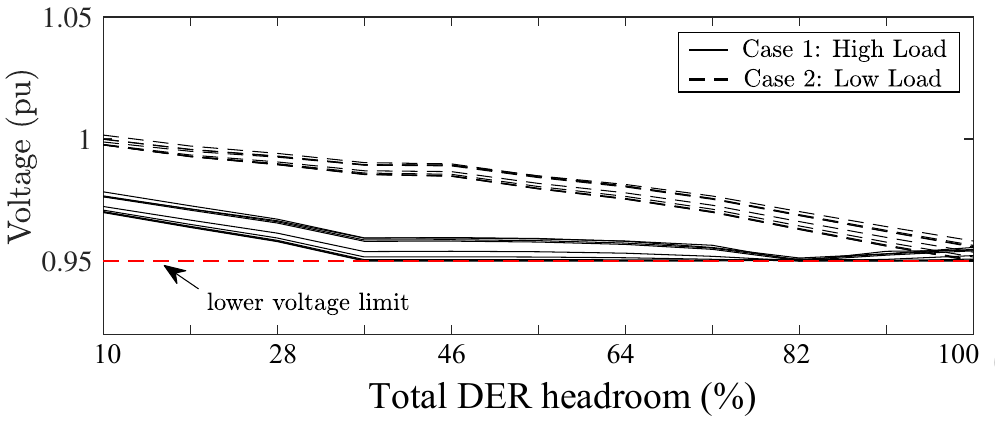}
\vspace{-2mm}
\caption {Voltage at all nodes downstream of node 733 in IEEE 37 bus system for case 1 and case 2 while estimating inductive capability curve ($\overline{q}_{net}^{cap}$) }
\label{fig:voltage_cap}
\end{figure}
Let's consider two cases with different loading conditions to compare the aggregated net capability curves i.e. high loading case 1 with peak load and low loading case 2 with half of the peak load. The capability curve for case 1 is shown in \figurename \ref{fig:cap_highload} as function of DER headroom with black solid lines. The shaded region enveloped by the black solid lines is the P-Q flexibility or capability domain; and, the dashed black line is base var demand, $q^{base}_{net}$ = 2.12 Mvar . The region above and below base var demand line can be seen as inductive and capacitive var support region respectively. Essentially, any point in the flexibility domain can be achieved by appropriate headroom. The RPFR ($[a, b]$) values for both case 1 and case 2 are computed using (\ref{eq:RPFR}) and are compared in Table \ref{tab:RPFR} for different DER real power headroom levels (e.g., for Case 1, 10\%: a = 0.34 - 2.12 = -1.78 and b = 4.01-2.12 = 1.89). It can be seen that the capacitive support region (magnitude of $a$) increases with increasing the headroom for both case 1 and 2 as increasing real power headroom frees the inverter capacity as well as it reduces the voltages due to increase in net load. This provides more scope for DERs to supply var leading to higher magnitude of $a$. However, the inductive var support region  (magnitude of $b$)  first increases with the headroom but starts decreasing towards the end for case 1 while for case 2, it continuously increases. This is because both the increasing the headroom and inductive var support cause low voltages and after a certain headroom level, the voltage of at least one node reaches to its minimum limit. Whereas, in case 2, the voltages do not reach to the minimum limit due to low load condition for almost all the headroom percentages as shown in \figurename \ref{fig:voltage_cap}. 


\begin{table}
\caption{RPFR (Mvar) at various DER headroom levels for Case 1 (peak load) and Case 2 (low load)}
\label{tab:RPFR}
\centering
\renewcommand{\arraystretch}{1.3}
\begin{tabular}{c  c  c  c  c} 
\hline
DER	& 10\% & 45\% & 65\% & 85\%  \\ 
\hline
Case 1 & [-1.78, 1.89]& [-3.48, 3.89]& [-3.84, 4.31]& [-4.04, 3.88]\\
Case 2 & [-1.72, 1.83]&[-3.36, 3.77]&[-3.71, 4.22]&[-3.91, 4.47]\\ 
\hline
\end{tabular}
\end{table}

\color{blue}

\begin{figure}
\color{blue}
\centering
\includegraphics[clip, width=0.85\columnwidth]{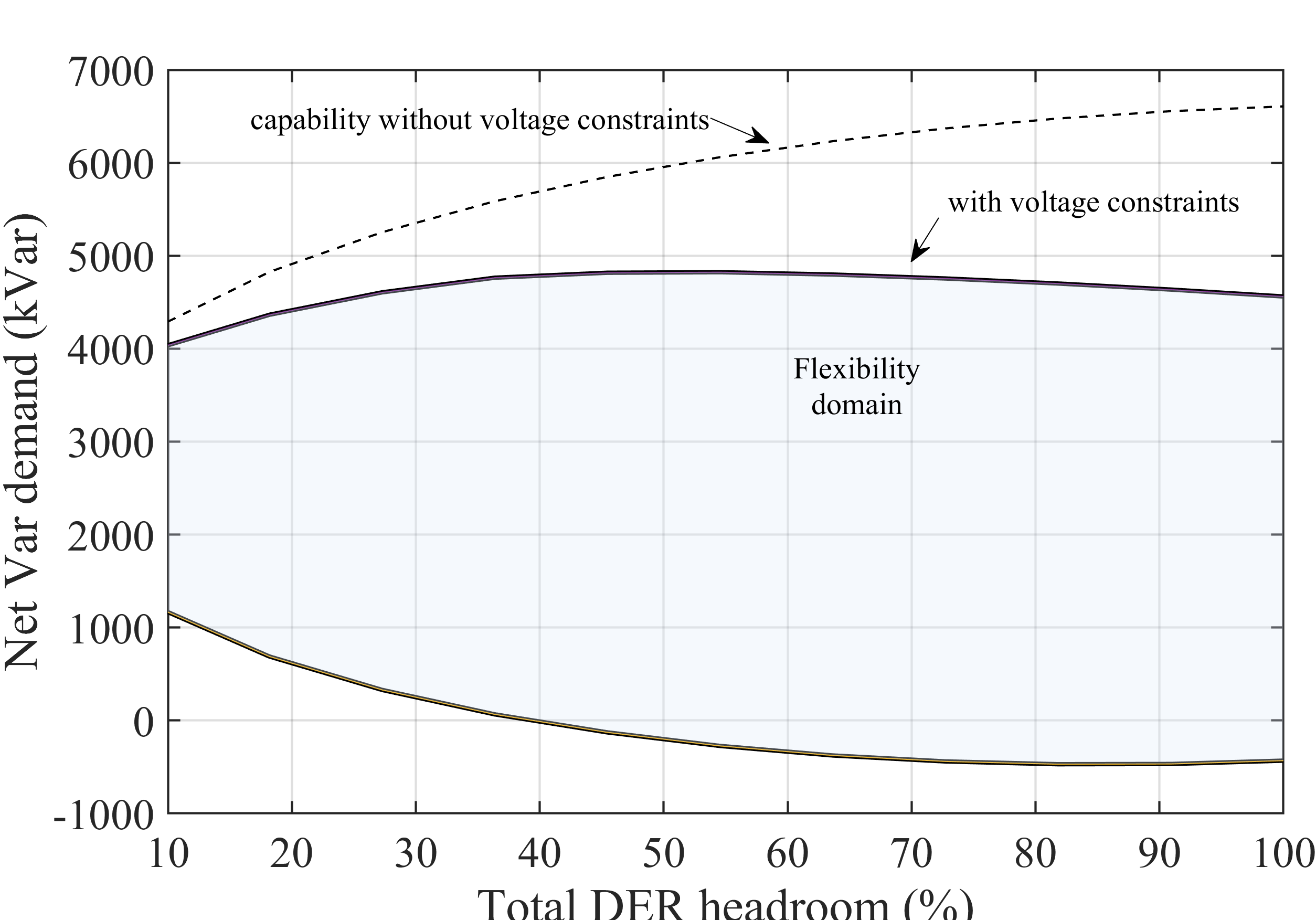}
\vspace{-2mm}
\caption {The aggregated net capability curve of IEEE 123 node DS with high DER penetration as function of DER headroom} 
\vspace{-0mm}
\label{fig:cap_123}
\end{figure}

Similar to 37 node system, the capability curve for IEEE 123 system is shown in \figurename \ref{fig:cap_123}. The dotted line shows the total inverter capacity to provide reactive power if voltage constraints are not considered. Whereas a proper consideration of voltage constraints in distribution system leads to a reduced capability region in 123 node system.

\subsection{Validation of Capability Curves}
\begin{figure}
    \color{blue}
	\centering
	\includegraphics[trim=0in 0in 0in 0in,width=0.9\columnwidth]{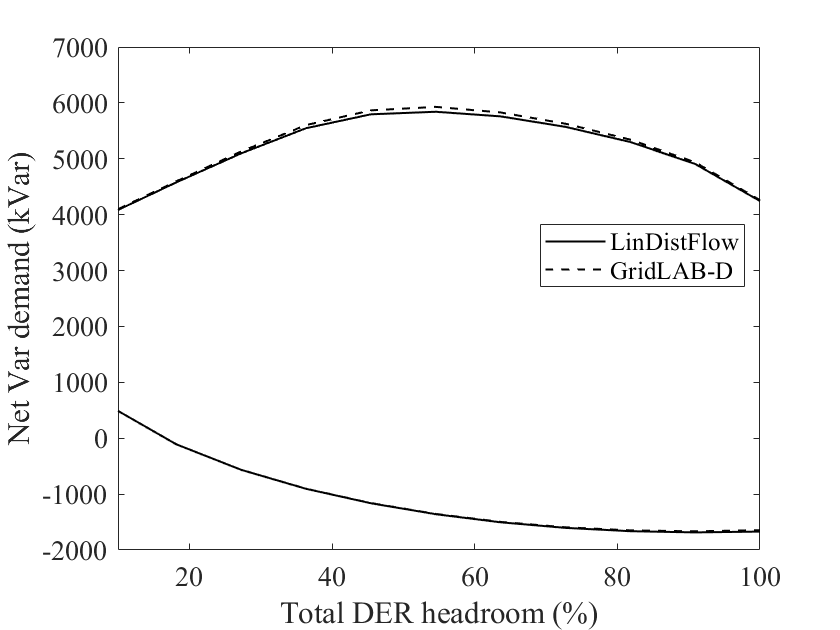}
    \caption {Validation of capability charts by comparing it against the full non-linear model (GridlabD) for IEEE 37 bus system}
    \label{fig:valid_cap}
\end{figure}

\begin{figure}
    \color{blue}
	\centering
	\includegraphics[trim=0in 0in 0in 0in,width=1.08\columnwidth]{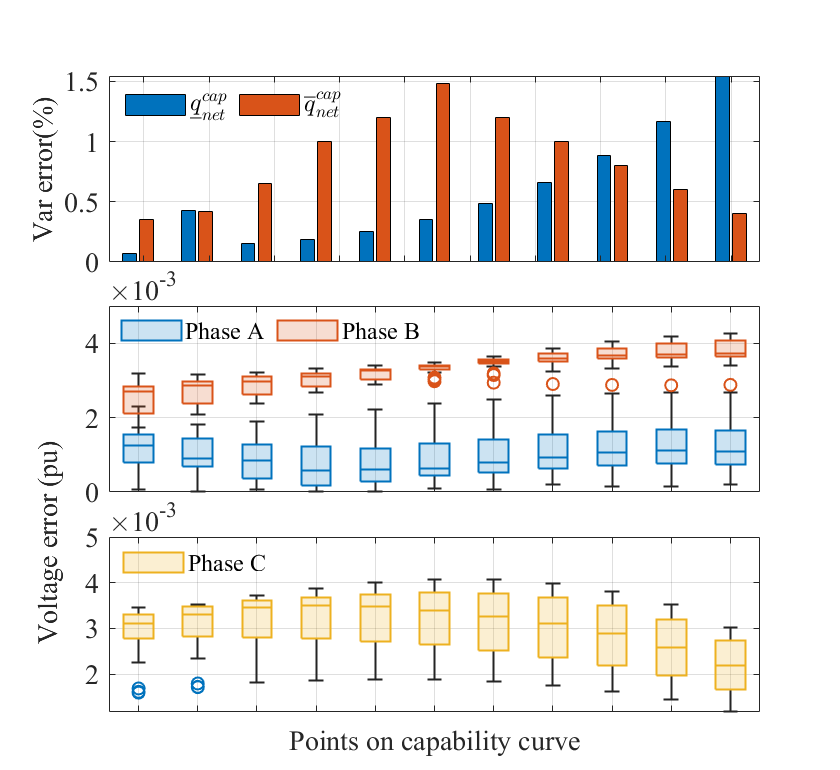}
    \caption {Validation of capability charts by comparing the DS voltages and substation var obtained from LinDistFlow optimization against the full non-linear model (GridlabD) for various points on the curve of IEEE 37 bus system}
    \label{fig:val}
\end{figure}

In order to validate the capability charts obtained from the LinDist3Flow model, we compare it against the full non-linear distribution system solver, GridLAB-D. \figurename \ref{fig:valid_cap} compares the capability curves obtained from the linearized model and GridLAB-D. The errors in var estimation for each point on both $\overline{q}^{cap}_{net}$ and $\underline{q}^{cap}_{net}$ curves are shown in \figurename \ref{fig:val} (top) for 37 bus system. The maximum error in var estimation is around 1.5\% whereas the average error is less than 1\% for both test systems as tabulated in Table \ref{tab:lin_dist_validation}. Further, it is important to check the accuracy of voltages while operating in the capability region. We check the voltage errors across all nodes and phases for various points on the upper and lower capability curves. The maximum and average errors in voltages are within 0.5\% and 0.3\% respectively. Note that the voltage errors are more than the base case since we are operating on extreme points on the capability curves. \figurename \ref{fig:val} (middle and bottom) show a box plot of distribution of errors in voltage across all nodes for each phase. It can be observed that the phase a has relatively smaller error in voltages.


\begin{table}[]
\color{blue}
\caption{Accuracy of estimated capability regions compared to non-linear solver  }
\label{tab:lin_dist_validation}
\centering
\renewcommand{\arraystretch}{1.2}
\begin{tabular}{cccccc}
\hline
\multirow{2}{*}{Test System} &  & \multicolumn{2}{c}{Voltage error (\%)} & \multicolumn{2}{c}{Substation var error   (\%)} \\ \cline{3-6} 
 &  & mean & max & mean & max \\ \hline
\multirow{2}{*}{IEEE 37} & $\underline{q}^{cap}_{net}$ & 0.24 & 0.4 & 0.56 & 1.54 \\
 & $\overline{q}^{cap}_{net}$ & 0.3 & 0.48 & 0.82 & 1.48 \\
\multirow{2}{*}{IEEE 123} & $\underline{q}^{cap}_{net}$ & 0.29 & 0.5 & 0.42 & 1.64 \\
 & $\overline{q}^{cap}_{net}$ & 0.35 & 0.62 & 0.92 & 1.71 \\ \hline
\end{tabular}
\end{table}

\subsection{Impact of Unbalanced DER on Capability Curves}
\begin{figure}
\color{blue}
\centering
\includegraphics[trim=-0in 0in 0in 0in,width=0.9\columnwidth]{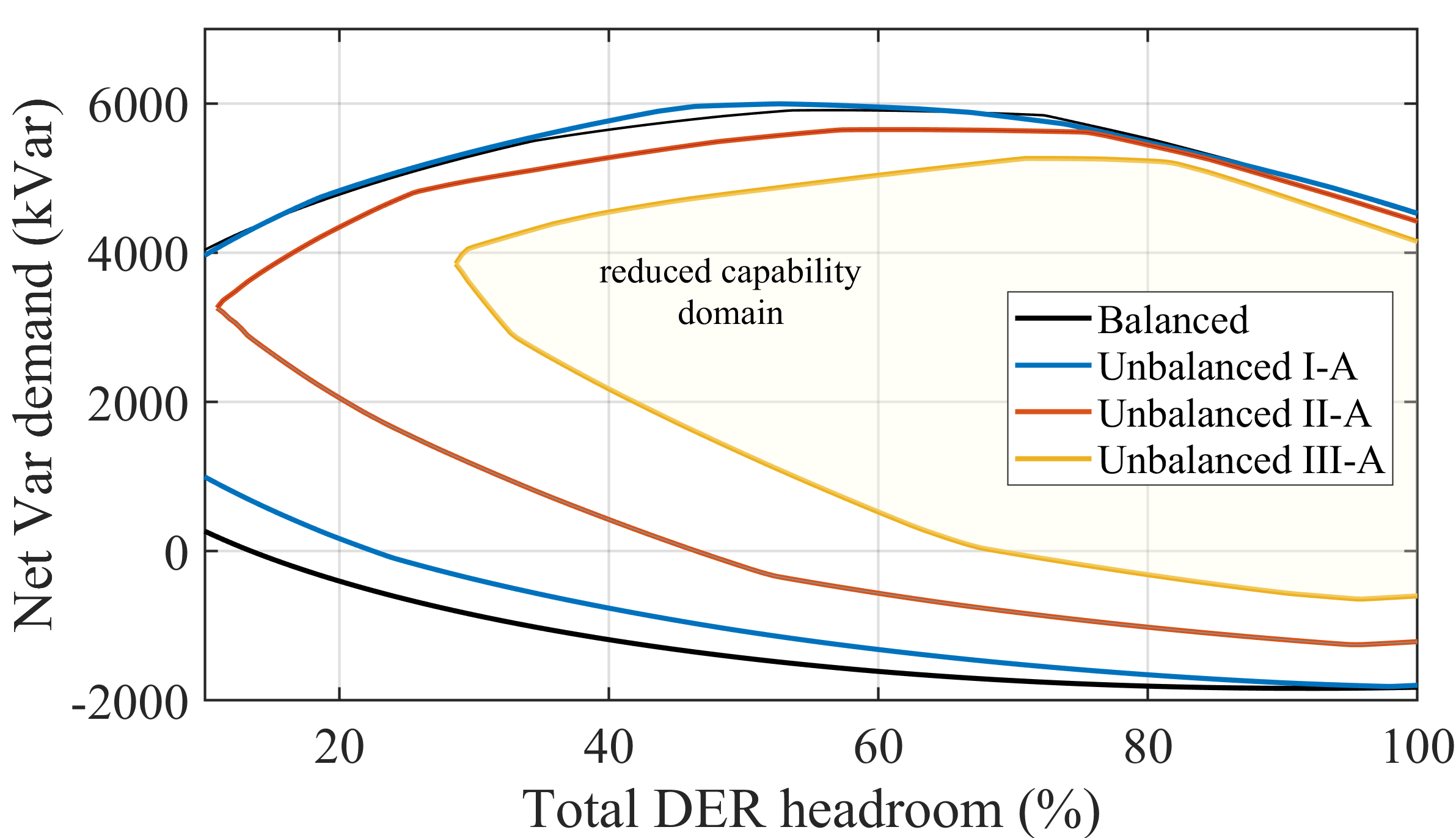}
\includegraphics[trim=-0in 0in 0in 0in,width=0.9\columnwidth]{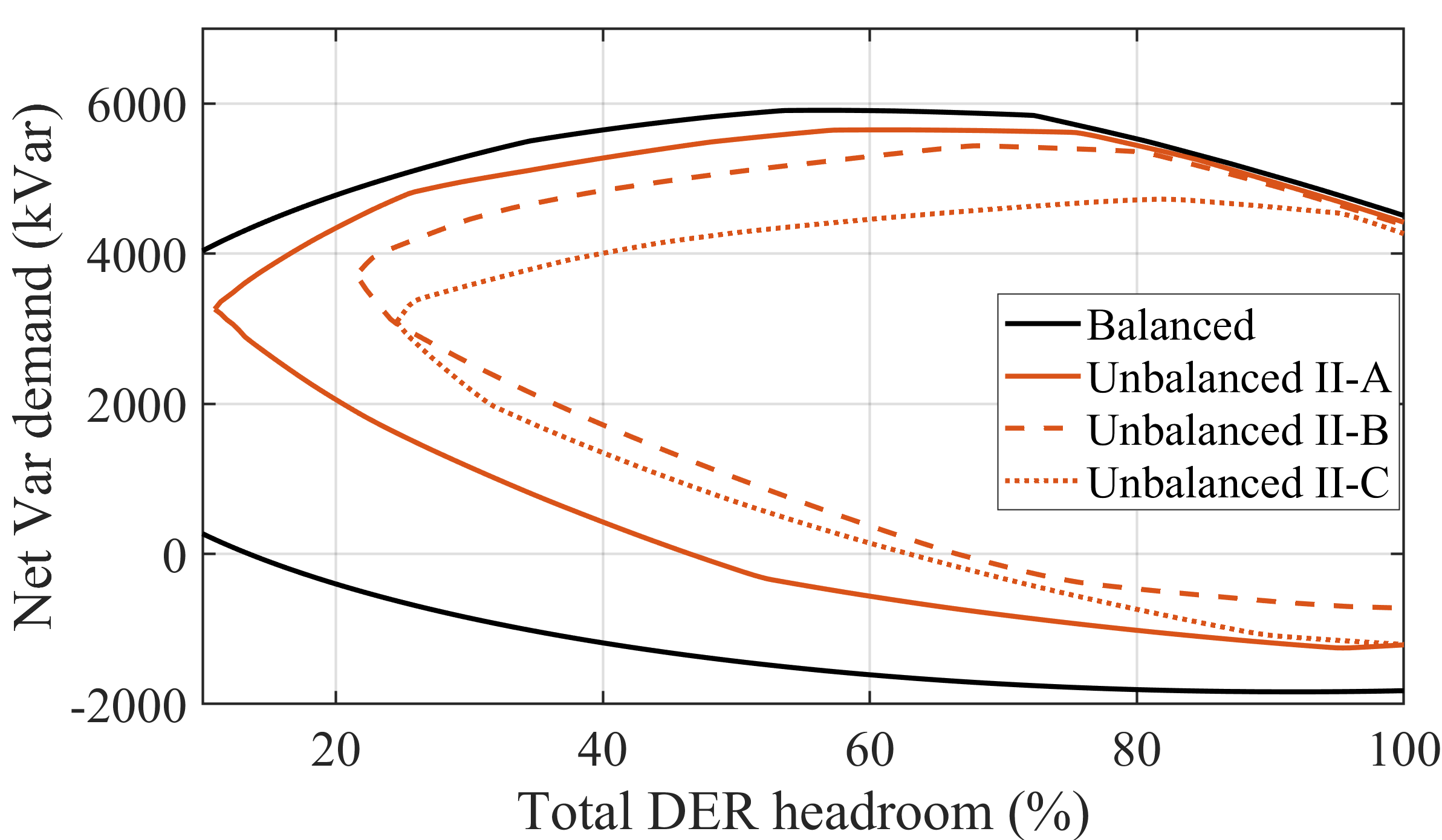}
\vspace{-2mm}
\caption {Impact of unbalance in DER distribution on the capability curve charts under various scenarios i.e. cases with (Top): increasing level of unbalance in phase A; (bottom): similar level of unbalance in different phases of the system}
\vspace{-0mm}
\label{fig:unbalanced}
\end{figure}

DER distribution in all 3 phases of a typical distribution system may be significantly unbalanced which in turn affect it's aggregated var capability. To demonstrate, we create various scenarios with increasing level of unbalance in DER distribution at all phases while keeping the total DER capacity (penetration level) same as shown in \ref{tab:unb}. Balanced scenario has equal distribution of DERs in all 3 phases. Unbalanced I, II, III scenario have increasing level of unbalance due to higher DER allocation at one of the phases. These scenario have been created for each phase one at a time.
\figurename \ref{fig:unbalanced}(top) shows the capability charts for increasing level of unbalance in phase A. It can be seen that the higher unbalance significantly shrinks the capability region. This happens because unbalanced distribution of DERs increases the voltage unbalance in the system making it more voltage constrained in one of the phases. Therefore, the full inverter capacity can not be utilized in order to maintain voltages within the operating bounds.  \figurename \ref{fig:unbalanced} (bottom) shows that even similar level of unbalance in different phases affects the capability curves differently owing to unbalanced nature (impedance and load) of distribution systems. This observation is in line with the discussion in \cite{alok2019}, where increasing unbalance in DER penetration is shown to impact the voltage stability margin of the system.
This emphasises the importance of capturing the unbalance in capability estimation process as proposed in this work.

\begin{table}[]
\color{blue}
\caption{Scenarios with increasing unbalanced in DER distribution in each of the 3 phases }
\label{tab:unb}
\centering
\renewcommand{\arraystretch}{1.2}
\begin{tabular}{ccccc}
\hline
\multirow{2}{*}{Scenario} & \multicolumn{3}{l}{\% Distribution of total DER capacity} \\
 & Phase A & Phase B & Phase C \\
 \hline
Balanced & 33.33 & 33.33 & 33.33 \\
Unbalanced I - A   & 50      & 25      & 25  \\
Unbalanced II - A  & 60      & 20      & 20   \\
Unbalanced III - A & 70      & 15      & 15    \\
Unbalanced I - B   & 25      & 50      & 25   \\
Unbalanced II - B  & 20      & 60      & 20   \\
Unbalanced III - B & 15      & 70      & 15    \\
Unbalanced I - C   & 25      & 25      & 50    \\
Unbalanced II - C  & 20      & 20      & 60    \\
Unbalanced III - C & 15      & 15      & 70    \\
\hline
\end{tabular}
\end{table}

\begin{figure}
\centering
\includegraphics[trim=-0in 0in 0in 0in,width=2.5in]{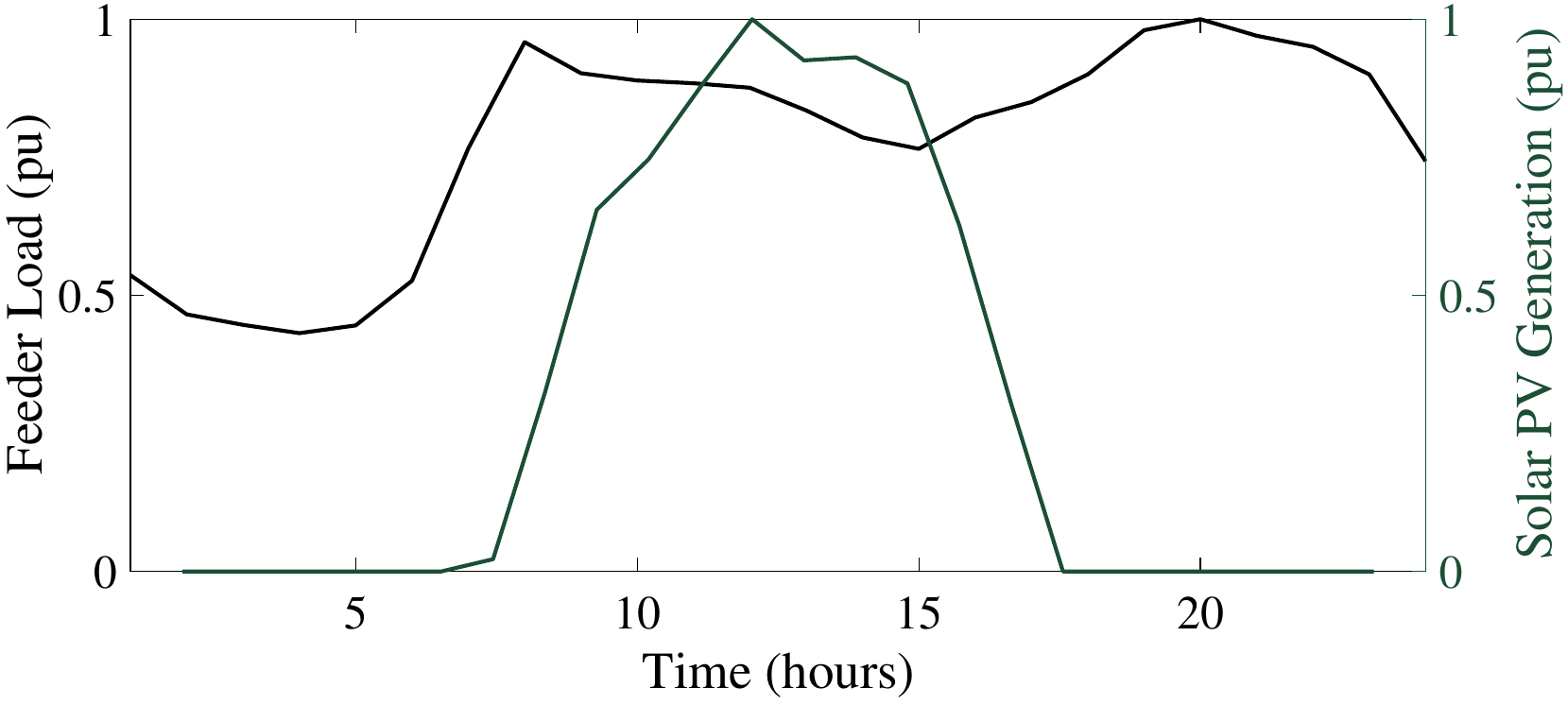}
\vspace{-3mm}
\caption {Normalized Daily load profile and solar PV generation profile for 24 hours with maximum value as 1 pu}
\label{fig:daily_profile}
\vspace{-4mm}
\end{figure}

\begin{figure}[b]
\centering
\vspace{-2mm}
\includegraphics[trim=-0in 0in 0in 0in,width=3.5in]{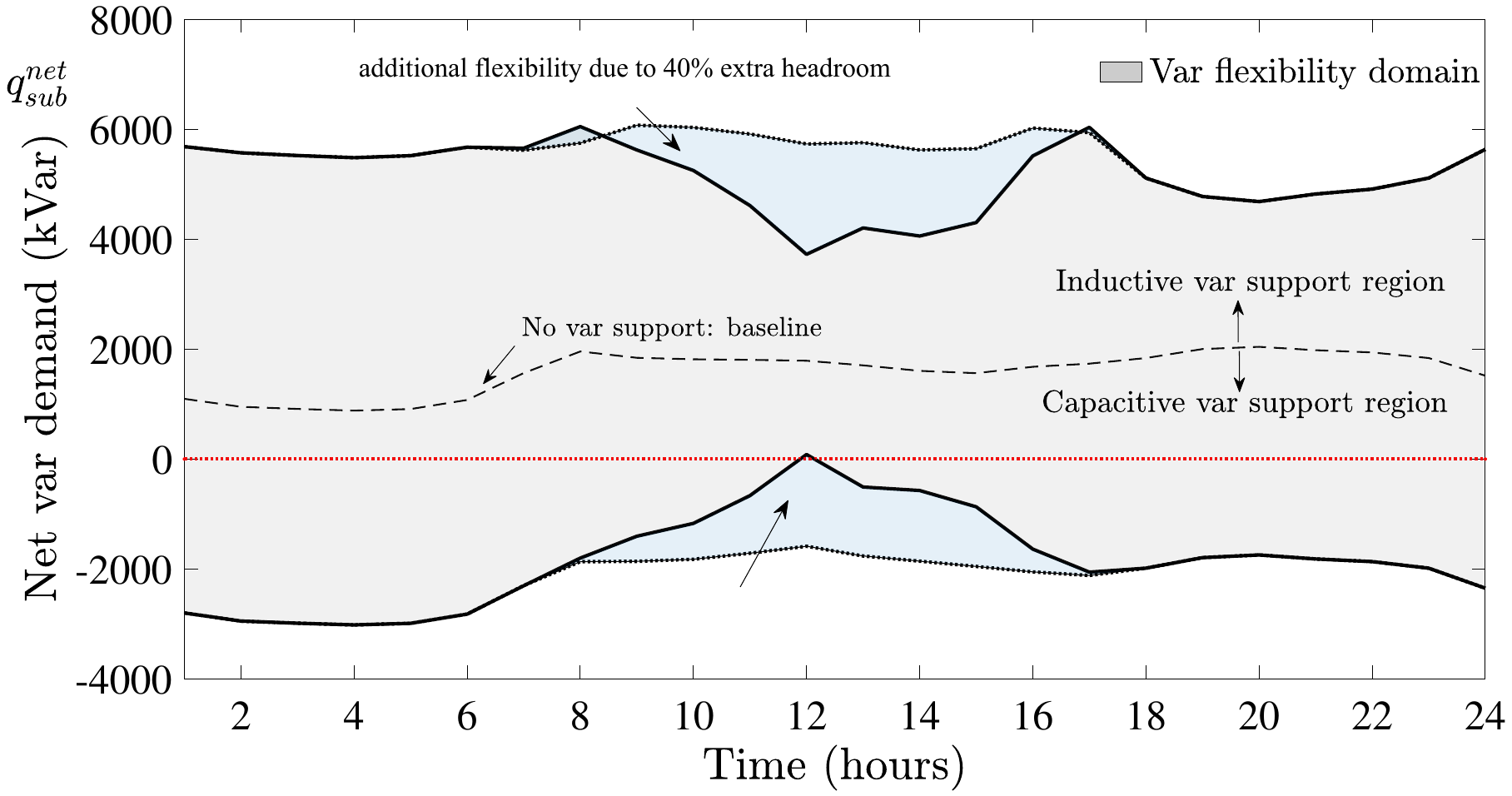}
\vspace{-5mm}
\caption {Day-ahead aggregated flexibility region of a DS with and without additional headroom}
\label{fig:day_cap}
\vspace{-0.5mm}
\end{figure}

\subsection{Impact of grid side voltages}
\begin{figure}
	\centering
	\color{blue}
    \vspace{-2mm}
	\includegraphics[trim=-0in 0in 0in 0in,width=3.5in]{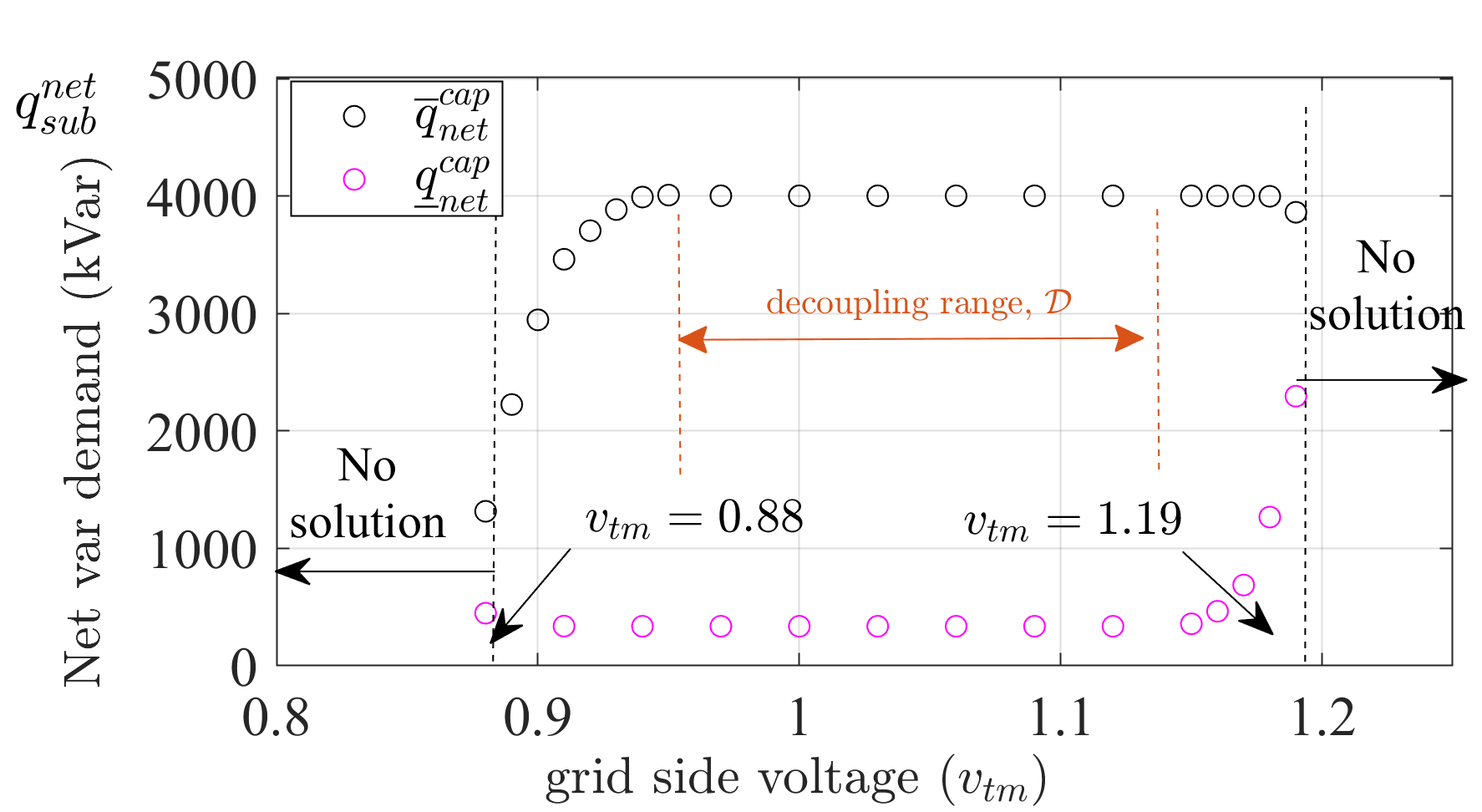}
    \vspace{-5mm}
    \caption {Impact of grid side voltage on the capability region}
    \label{fig:v_tm_impact}
\end{figure}

From transmission side, the primary substation voltage is a crucial factor that can affect the capability domain significantly. Though, usually we expect the $v_{tm}$ to be around 1, in case of contingencies and other events, it can significantly deviate from nominal value. \figurename \ref{fig:v_tm_impact} shows how the capability region varies with change in $v_{tm}$ at no curtailment. Note that the flexibility region shrinks as $v_{tm}$ moves away from nominal 1 pu on either side beyond the decoupling range $\mathcal{D}$. It can be seen that the DER-OPF becomes infeasible for $v_{tm}$ greater than 1.19 pu and less than 0.88 pu that means no flexibility is available without violating the voltage limits.


\color{black}
\subsection{Day-ahead Capability Curve}
In the last section, capability curve were shown for a given operating condition,
however, the day-ahead capability curves can also be estimated to be utilized  by TSO for day-ahead planning. 
A normalized daily load curve and solar PV generation profile is applied to each load and PV unit respectively as shown in \figurename \ref{fig:daily_profile}.  \figurename \ref{fig:day_cap} shows a capability curves (black solid lines) and var support region (grey shaded area) of the test system at hourly operating points with no headroom. The day-ahead curve gives more visual information of how the aggregated capability varies with changing operating condition throughout the day. {\color{blue}It can be seen that the flexibility range is minimum at noon when least inverter capacity is available for var support (between 0800 hrs to 1700 hrs). However, a 40\% headroom will free the inverter capacity and expand its flexibility area shown by the blue shaded portion during peak solar hours (0800 hrs to 1700 hrs) as shown in \figurename \ref{fig:day_cap}.}

\subsection{Factors Affecting Capability Curve}
\subsubsection{Impact of DER penetration levels}
Table \ref{tab:RPFR_pen} shows the $RPFR$ values for increasing DER penetration level with no headroom for the high load condition. As expected, both nominal capacitive and induction flexibility region increase with higher DER penetration. 

\begin{table}[b]
\vspace{-4mm}
\caption{RPFR (Mvar) for different DER penetration levels}
\label{tab:RPFR_pen}
\centering
\renewcommand{\arraystretch}{1.3}
\begin{tabular}{ c |  c  c  c  c  c  }
\hline
\multicolumn{1}{c}{DER penetration $\rightarrow$}	& 20\% & 40\% & 60\% & 80\% & 100\%  \\
\hline
$[\,a, $              & [-043, &[-0.85,&[-1.26,&[-1.67,&[-2.07,\\
$ b\,]$   &  0.4]& 0.87]& 1.31]&  1.76]&  2.22] \\
\hline
\end{tabular}
\end{table}

\subsubsection{Impact of Inverter size}
Inverter size plays a crucial role in available DER capability. Table \ref{tab:RPFR_invsize} compares the $RPFR$ for different inverter sizes during peak solar generation. Inverters with no oversize at peak generation resuls in no headroom leading to zero flexibility region. However, the headroom can be increased by freeing up the inverter capacity via either cutailment or storage as shown in the Table \ref{tab:RPFR_invsize} (row 2). To realize headroom, Table \ref{tab:RPFR_invsize} shows the trade-off between oversizing and curtailment, to achieve desired flexibility during peak solar generation. This trade-off is also relevant to comply to the integration standard 1547-2018 as discussed in the next section.  

\begin{table}
\caption{RPFR (Mvar) for different inverter sizing}
\label{tab:RPFR_invsize}
\centering
\renewcommand{\arraystretch}{1.3}
\begin{tabular}{ c  c  c  c   }
\hline
\multicolumn{1}{c}{Inverter oversize $\rightarrow$}	& 1 & 1.1 & 1.2 \\
\hline
No DER headroom                & [0, 0] & [-1.78, 1.89] & [-2.54, 2.76] \\
40\% DER headroom &  [-3.04, 3.35] & [-3.48, 3.89] & [-3.89, 4.42] \\
\hline
\end{tabular}
\end{table}

\subsection{Integration Standard IEEE1547 Compliance}
\begin{figure}
\centering
\vspace{-2mm}
\includegraphics[trim=-0in 0in 0in 0in,width=3.5in]{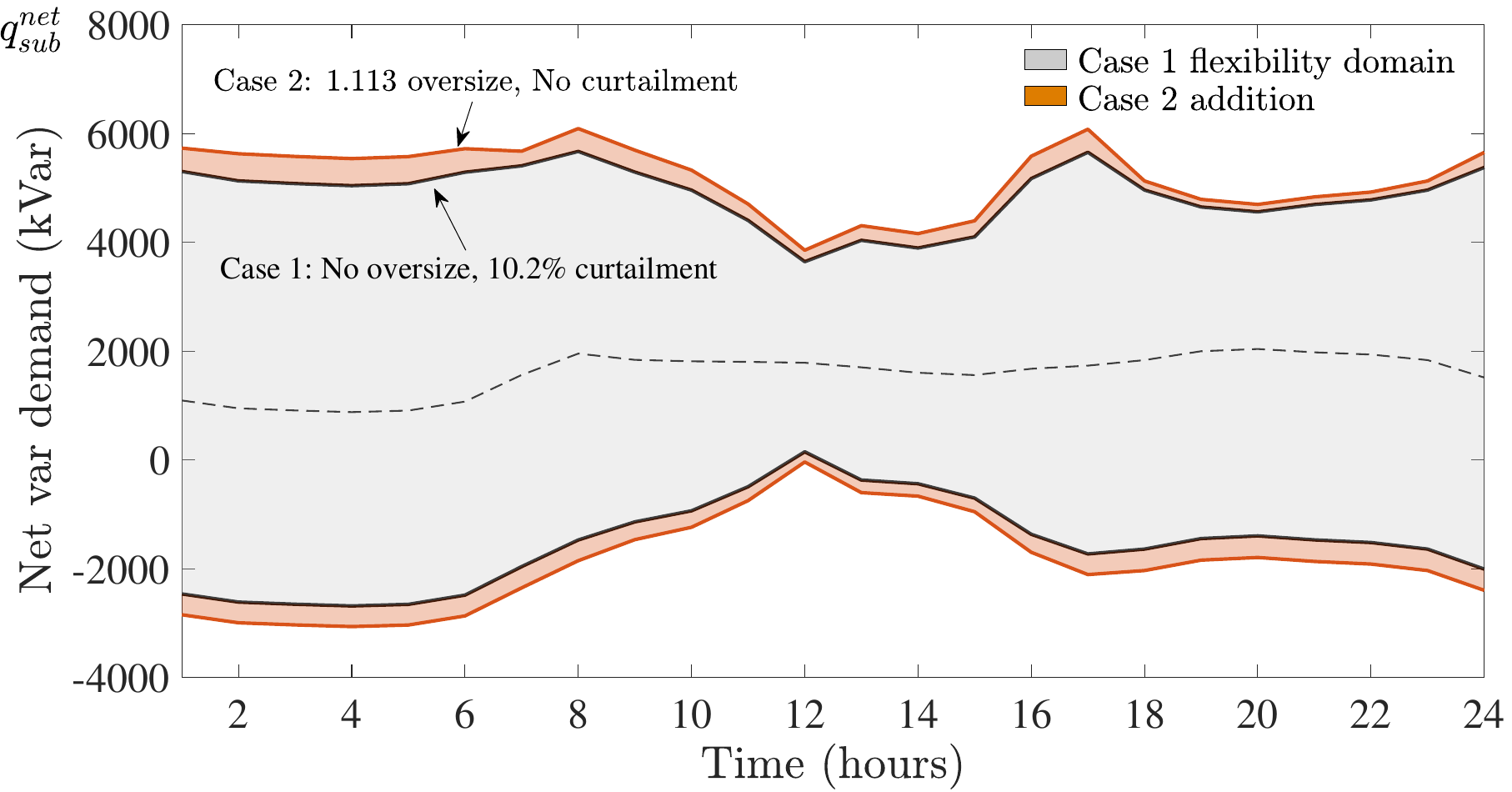}
\vspace{-5mm}
\caption {Two different day-ahead aggregated capability curves and domains due to compliance to IEEE1547 var capability requirements}
\vspace{-2mm}
\label{fig:1547}
\end{figure}

The recently revised DER integration standard IEEE1547-2018 has made it compulsory for each inverter-based DER unit to provide var capability of 44\% of its kW rating at all operating conditions \cite{noauthor_ieee_2018}. In order to comply with it, there are two possible options i.e. either oversize the inverter by 1.113 times the kW rating with no DER headroom or 10.2\% DER headroom during peak hours with no inverter oversize. Here, we have compared the impact of IEEE1547 compliance on the aggregated RPFR by choosing both the options in form of two cases as shown in \figurename \ref{fig:1547}. Case 1 is shown by black solid lines where 10.2\% headroom at each DER with no oversized inverters and case 2 is shown by orange solid lines  where each inverter is oversize by 1.113 times. Note that both the cases comply to IEEE1547 standard, however the case 2 has broader aggregated RPFR ([-1.89, 2.03]) compare to case 1 ( [-1.71, 1.82]) at 12 noon. At the same time, the case 1 might be more viable as the real power headroom is provided occasionally when the var is needed by TSO, whereas the inverter oversize cost is permanent. Nonetheless, it depends on many other factors such as policy, incentive structure, ancillary service market etc. and further cost-benefit analysis is needed of specific cases to arrive at any decision.

\section{Impact of Aggregated DER Var Support on the Transmission Grid}
\subsection{Impact of var Capability on Transmission System Generation}
Let us consider the IEEE 9-bus transmission system with the distributed inverters connected at load bus 7 (T7) (100 MW, 35 Mvar). We consider different scenarios for the IEEE 9 bus system with and without the var support available at the transmission system. We have considered 50\% smart inverter penetration with 10\% headroom (minimum headroom) and 44\% of its rating available for reactive power modulation for the cases with var support. The resultant power injection from the inverters at load bus 7 is (45 MW, 22 Mvar) (with var support). The case without var support assumes the inverters operate in unity power factor (UPF) mode, which means there is no reactive power injection from inverters. The 
The individual var generation from the three synchronous generators in IEEE 9- bus system are shown in Table IV. In this system, for the contingency of line 5-6 outage, we demonstrate that with the var capability available to the transmission system, the total reactive power generation is reduced significantly. In cases of no contingency also, the case for No var support has the higher reactive power generation. 

\begin{table}
\caption{DER var Impact of Transmission var Generation}
\label{tab:ts-Qgen}
\centering
\renewcommand{\arraystretch}{1.3}
\begin{tabular}{ c  c  c  c   }
\hline
\multicolumn{1}{c}{Scenario } 	& $Q_{Gen1}$ & $Q_{Gen2}$ & $Q_{Gen3}$ \\
                     &(Mvar) & (Mvar) & (Mvar)\\
\hline
50\% DER (No var, No Contingency)               & 30.28 & 14.92 & -3.52 \\
50\% DER (With var, No Contingency) &  26.4 & 3.75 & -12.95 \\
50\% DER (No var, Contingency) &  85.62 & 45.36 & 10.21 \\


50\% DER (With var, Contingency) &  82.30 & 33.37 & -0.31 \\
\hline
\end{tabular}
\end{table}

From Table \ref{tab:ts-Qgen} , for each case (contingency or no contingency), the total reactive power generation from the conventional generators is lower for the cases with the inverter var capability available to the TSO. For the case with the inverters just supplying real power, the reactive power generation is the highest under the contingency case. The next subsection discusses the impact of var capability of the DSO's by modeling the T-D systems in detail using T-D co-simulation methodology \cite{sun_master_2015}- \cite{alok_cosim_2021}.

{
\color{blue}The present framework discusses the availability of another var resource through DSO flexibility. However an economic analysis with the reactive pricing can provide further inputs that enables the TSO to make optimal choice. In most disturbances, availability of local var resources through DSO capability can prove more effective compared to central var dispatch from large synchronous generators.}

\subsection{Utilizing var Capability through TSO-DSO Interactions}
\begin{figure}
	\centering
    \vspace{-2mm}
	\includegraphics[trim=-0in 0in 0in 0in,width=2in]{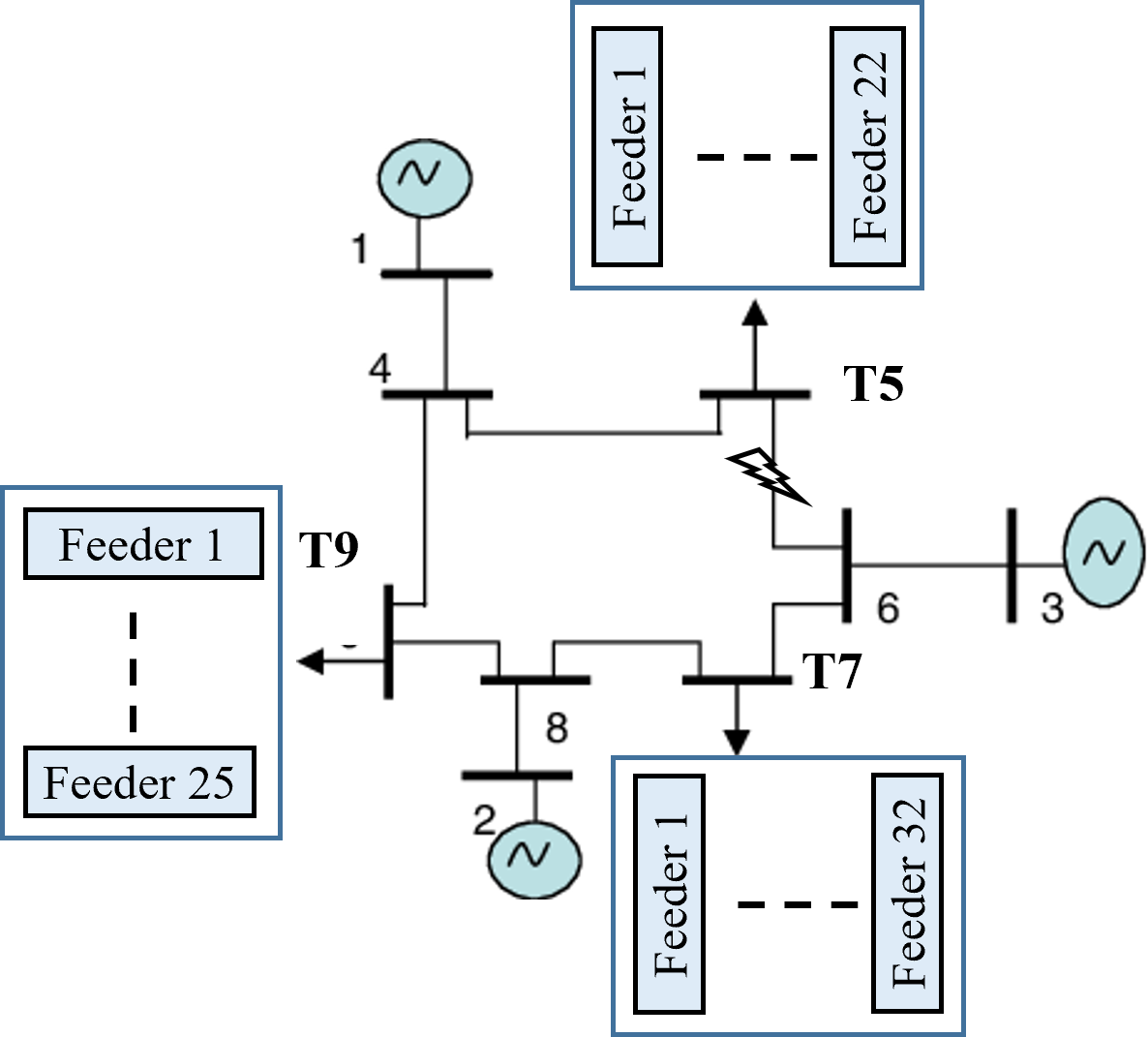}
    \vspace{-4mm}
    \caption {Integrated T-D test system with coupled IEEE 9 bus TS and multiple IEEE 37 bus DS feeders}
    \label{fig:T-DSystem}
\end{figure}

The proposed framework estimates the aggregated var capability curve for the transmission grid. However, the grid might not need the maximum var support all the time; rather it can ask for the var support in specific needs e.g. in case of voltage dips due to line contingencies. In this section, we will demonstrate on an integrated T-D test system, how the proposed aggregated var capability can potentially enhance the options for TSO. A T-D cosimulation platform is developed based on reference \cite{sun_master_2015, alok_cosim_2021} to accurately model the T-D interactions. An integrated T-D test system is modeled by coupling aggregated multiple  IEEE 37 bus DS feeders (22,32,25) at three load buses (T5,T7,T9) of the IEEE 9 bus TS as shown in \figurename \ref{fig:T-DSystem}. 

\begin{figure}
	\centering
	\includegraphics[trim=-0in 0in 0in 0in,width=3.5in]{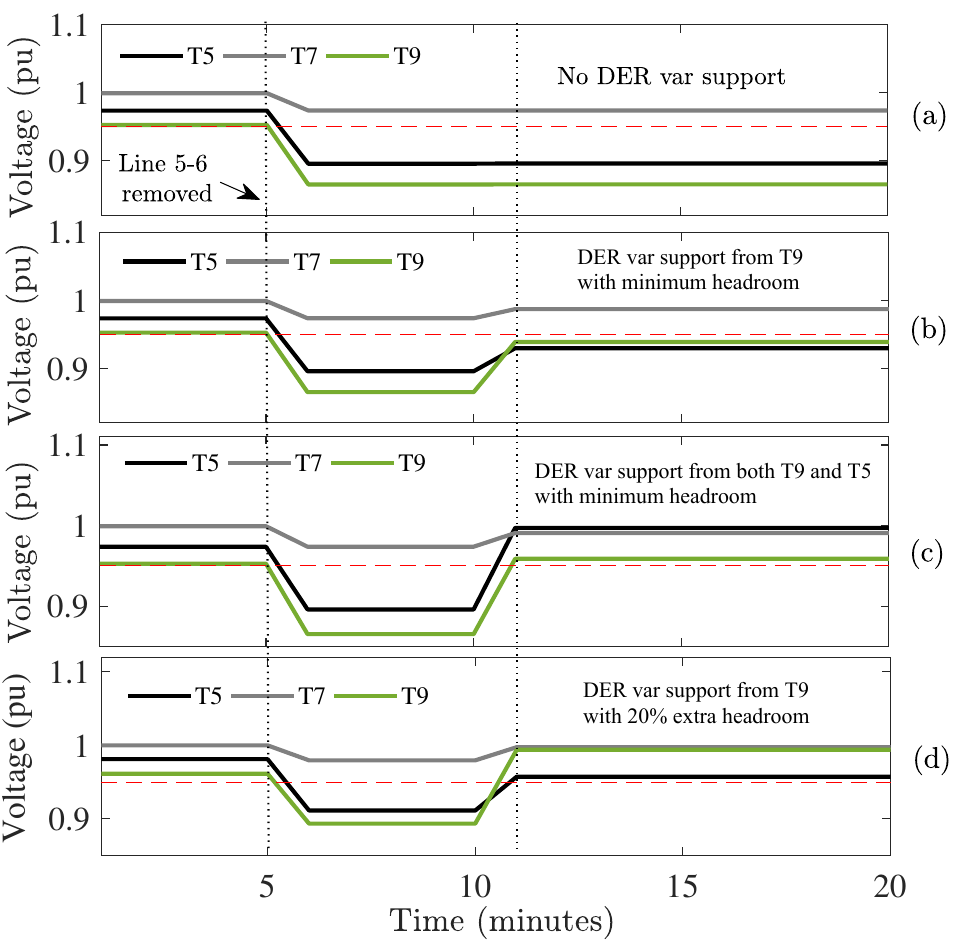}
    \vspace{-3mm}
    \caption {Impact of Aggregated DER var support on grid voltages after transmission line 5-6 contingency in various cases}
    \label{fig:contingency56_1}
    \vspace{-3mm}
\end{figure}

Let's consider a operating point with peak solar generation to demonstrate the impact of minimum available var flexibility. We will compare the impact of DER var support under line T5-6 contingency for following cases: a) No DER var support provided by any DSO; b) DSO at bus T9 provides just enough DER var support to comply with 1547; c) DSO at both bus T9 and T5 provide just enough DER var support to comply with integration standard 1547; d) DSO at bus T9 provide more DER var support than case (b) by 20\% headroom. 
\figurename \ref{fig:contingency56_1} compares the voltages at transmission buses for all cases.  At $t=5$, line 5-6 is removed that leads to dip in voltages and bus T5 and T9 suffer under voltage violation. In case (b), the support by only T9 is not enough to recover voltages above 0.95. In such cases, TSO either can request var flexibility from both T9 and T5 as recommended by 1547 standard i.e. case (c) or it can request extra support from T9 that can be provided by some headroom i.e., case (d). It can be seen that both case (c) and (d) recover voltages above the limit, however, the amount of voltage boost at T9 and T5 differ based on the cases. 

From the results shown in \figurename \ref{fig:contingency56_1}, we can see that in case of no DERs, under a contingency, the voltages at the load buses are violated with respect the lower limit of 0.95 pu. With minimum headroom to comply to the IEEE 1547-2018, var capability at T9 is used, but the voltages are not able to recover to more than 0.95 pu. If the DSO does not have provisions to change the headroom, then the neighboring DSO can utilize its var capability to help the TSO as shown in results in (c). However, if the DSO can adjust the DER headroom, then the DSO at T9 can increase the headroom to further enhance the amount of var support to TSO as shown in (d).

Note that the estimation of optimal var support request profile depends on the various factor such as objective of TSO, availability of DER flexibility, economic compensation policies etc and needs to be achieved via an optimization process which is beyond the scope of this paper. The 4 cases here demonstrate the potential of the proposed framework that provides higher flexibility to TSO. {\color{blue}The var capability curve can be used as an input by the TSO to optimize its operation by varying the var within the var capability curve. The var capability curve allows for TSO to request a DSO to dispatch a var at the TSO-DSO interface}

\vspace{-5mm}


\section{Conclusion}
This work is one of the components in the pursuit of utilizing the increasing DER penetration for the benefit of the future grid, driven by the proposition that multitude of geographically distributed DERs with var control capability can be seen as flexible var resources (mini SVCs) for the grid. To facilitate this vision, a systematic OPF based methodology is proposed to construct an \textit{aggregated net var capability curve} of a DS with high DER penetration, analogous to a conventional bulk generator. The proposed capability curve also accounts for DER headroom that enables TSO to utilize both P and Q flexibility provided by DERs into their planning and operational activities. 
The impact of DS voltage constraints, inverter sizing and, T-D coupling on the flexibility region are discussed via results on an unbalanced IEEE 37 bus DS.  
Finally, It is shown how aggregated DER var flexibility can affect the transmission system performance on an integrated T-D test system in a cosimulation environment. 

Certainly, the formulation and details of inclusion of the proposed flexibility in grid planning and operations remain an exciting challenge for future studies that require a larger discussion on policy, payment structure etc. Nonetheless, the results are encouraging and indicate that the proposed aggregated capability indeed has potential to improve grid optimality by providing enhanced flexibility services to TSO.

\vspace{-2mm}
\bibliographystyle{IEEEtran}
\bibliography{IEEEabrv,References_zotero}

\begin{thebibliography}{10}
\providecommand{\url}[1]{#1}
\csname url@samestyle\endcsname
\providecommand{\newblock}{\relax}
\providecommand{\bibinfo}[2]{#2}
\providecommand{\BIBentrySTDinterwordspacing}{\spaceskip=0pt\relax}
\providecommand{\BIBentryALTinterwordstretchfactor}{4}
\providecommand{\BIBentryALTinterwordspacing}{\spaceskip=\fontdimen2\font plus
\BIBentryALTinterwordstretchfactor\fontdimen3\font minus
  \fontdimen4\font\relax}
\providecommand{\BIBforeignlanguage}[2]{{%
\expandafter\ifx\csname l@#1\endcsname\relax
\typeout{** WARNING: IEEEtran.bst: No hyphenation pattern has been}%
\typeout{** loaded for the language `#1'. Using the pattern for}%
\typeout{** the default language instead.}%
\else
\language=\csname l@#1\endcsname
\fi
#2}}
\providecommand{\BIBdecl}{\relax}
\BIBdecl

\bibitem{bao_online_2003}
L.~Bao, Z.~Huang, and W.~Xu, ``Online voltage stability monitoring using {VAr}
  reserves,'' \emph{IEEE Transactions on Power Systems}, vol.~18, no.~4, pp.
  1461--1469, Nov. 2003.

\bibitem{song_reactive_2003}
H.~Song, B.~Lee, S.-H. Kwon, and V.~Ajjarapu, ``Reactive reserve-based
  contingency constrained optimal power flow ({RCCOPF}) for enhancement of
  voltage stability margins,'' \emph{IEEE Transactions on Power Systems},
  vol.~18, no.~4, pp. 1538--1546, Nov. 2003.

\bibitem{goergens_determination_2015}
P.~Goergens, F.~Potratz, M.~Gödde, and A.~Schnettler, ``Determination of the
  potential to provide reactive power from distribution grids to the
  transmission grid using optimal power flow,'' in \emph{{International}
  {Universities} {Power} {Engineering} {Conference}}, Sep. 2015, pp. 1--6.

\bibitem{barth_technical_2013}
H.~Barth, D.~Hidalgo, A.~Pohlemann, M.~Braun, L.~H. Hansen, and H.~Knudsen,
  ``Technical and economical assessment of reactive power provision from
  distributed generators: {Case} study area of {East} {Denmark},'' in
  \emph{2013 {IEEE} {Grenoble} {Conference}}, Jun. 2013, pp. 1--6.

\bibitem{keane_state---art_2013}
A.~Keane, L.~F. Ochoa, C.~L.~T. Borges, G.~W. Ault, A.~D. Alarcon-Rodriguez,
  R.~A.~F. Currie, F.~Pilo, C.~Dent, and G.~P. Harrison, ``State-of-the-{Art}
  {Techniques} and {Challenges} {Ahead} for {Distributed} {Generation}
  {Planning} and {Optimization},'' \emph{IEEE Transactions on Power Systems},
  vol.~28, no.~2, pp. 1493--1502, May 2013.

\bibitem{perez-arriaga_transmission_2016}
I.~J. Perez-Arriaga, ``The {Transmission} of the {Future}: {The} {Impact} of
  {Distributed} {Energy} {Resources} on the {Network},'' \emph{IEEE Power and
  Energy Magazine}, vol.~14, no.~4, pp. 41--53, Jul. 2016.

\bibitem{singhal_real-time_2018}
A.~Singhal, V.~Ajjarapu, J.~C. Fuller, and J.~Hansen, ``Real-{Time} {Local}
  {Volt}/{VAR} {Control} {Under} {External} {Disturbances} with {High} {PV}
  {Penetration},'' \emph{IEEE Transactions on Smart Grid}, pp. 1--1, 2018.

\bibitem{zhu_fast_2016}
H.~Zhu and H.~J. Liu, ``Fast {Local} {Voltage} {Control} {Under} {Limited}
  {Reactive} {Power}: {Optimality} and {Stability} {Analysis},'' \emph{IEEE
  Transactions on Power Systems}, vol.~31, no.~5, pp. 3794--3803, Sep. 2016.

\bibitem{zhang_optimal_2015}
B.~Zhang \emph{et~al.}, ``An {Optimal} and {Distributed} {Method} for {Voltage}
  {Regulation} in {Power} {Distribution} {Systems},'' \emph{IEEE Trans. on
  Power Sys.}, vol.~30, no.~4, pp. 1714--1726, Jul. 2015.

\bibitem{noauthor_ieee_2018}
``{IEEE} {Standard} for {Interconnection} and {Interoperability} of
  {Distributed} {Energy} {Resources} with {Associated} {Electric} {Power}
  {Systems} {Interfaces},'' \emph{IEEE Std 1547-2018}, pp. 1--138, Apr. 2018.

\bibitem{noauthor_impact_2018}
``\BIBforeignlanguage{en-US}{Impact of {IEEE} 1547 {Standard} on {Smart}
  {Inverters}},'' IEEE PES Industry Technical Support Task Force, Technical
  {Report} PES-TR67, May 2018.

\bibitem{marten_analysis_2013}
F.~Marten \emph{et~al.}, ``Analysis of a reactive power exchange between
  distribution and transmission grids,'' in \emph{2013 {IEEE} {International}
  {Workshop} on {Inteligent} {Energy} {Systems} ({IWIES})}, Nov. 2013, pp.
  52--57.

\bibitem{konopinski_extended_2009}
R.~J. Konopinski, P.~Vijayan, and V.~Ajjarapu, ``Extended {Reactive}
  {Capability} of {DFIG} {Wind} {Parks} for {Enhanced} {System}
  {Performance},'' \emph{IEEE Trans. on Power Sys.}, vol.~24, no.~3, pp.
  1346--1355, Aug. 2009.

\bibitem{cuffe_transmission_2012}
P.~Cuffe, P.~Smith, and A.~Keane, ``Transmission {System} {Impact} of {Wind}
  {Energy} {Harvesting} {Networks},'' \emph{IEEE Transactions on Sustainable
  Energy}, vol.~3, no.~4, pp. 643--651, Oct. 2012.

\bibitem{kundu_approximating_2018}
S.~Kundu, K.~Kalsi, and S.~Backhaus, ``Approximating {Flexibility} in
  {Distributed} {Energy} {Resources}: {A} {Geometric} {Approach},'' in
  \emph{2018 {Power} {Systems} {Computation} {Conference} ({PSCC})}, Jun. 2018,
  pp. 1--7.

\bibitem{marten_optimizing_2014}
F.~Marten \emph{et~al.}, ``Optimizing the reactive power balance between a
  distribution and transmission grid through iteratively updated grid
  equivalents,'' in \emph{{Power} {Systems} {Computation} {Conf.}}, Aug. 2014,
  pp. 1--7.

\bibitem{kaempf_reactive_2014}
E.~Kaempf \emph{et~al.}, ``Reactive power provision by distribution system
  operators — {Optimizing} use of available flexibility,'' in \emph{{IEEE}
  {PES} {Innovative} {Smart} {Grid} {Technologies}, {Europe}}, Oct. 2014, pp.
  1--5.

\bibitem{8442917}
F.~Capitanescu, ``Ac opf-based methodology for exploiting flexibility provision
  at tso/dso interface via oltc-controlled demand reduction,'' in \emph{2018
  Power Systems Computation Conference (PSCC)}, 2018, pp. 1--6.

\bibitem{8291006}
J.~Silva, J.~Sumaili, R.~J. Bessa, L.~Seca, M.~A. Matos, V.~Miranda,
  M.~Caujolle, B.~Goncer, and M.~Sebastian-Viana, ``Estimating the active and
  reactive power flexibility area at the tso-dso interface,'' \emph{IEEE
  Transactions on Power Systems}, vol.~33, no.~5, pp. 4741--4750, 2018.

\bibitem{9543347}
G.~C. Kryonidis, A.~N. Lois, K.-N.~D. Malamaki, and C.~S. Demoulias,
  ``Two-stage approach for the provision of time-dependent flexibility at
  tso-dso interface,'' in \emph{2021 International Conference on Smart Energy
  Systems and Technologies (SEST)}, 2021, pp. 1--6.

\bibitem{CAPITANESCU2018226}
\BIBentryALTinterwordspacing
F.~Capitanescu, ``Tso–dso interaction: Active distribution network power
  chart for tso ancillary services provision,'' \emph{Electric Power Systems
  Research}, vol. 163, pp. 226--230, 2018. [Online]. Available:
  \url{https://www.sciencedirect.com/science/article/pii/S0378779618301822}
\BIBentrySTDinterwordspacing

\bibitem{9295337}
B.~Cui, A.~Zamzam, and A.~Bernstein, ``Network-cognizant time-coupled aggregate
  flexibility of distribution systems under uncertainties,'' \emph{IEEE Control
  Systems Letters}, vol.~5, no.~5, pp. 1723--1728, 2021.

\bibitem{Kara}
S.~Karagiannopoulos, C.~Mylonas, P.~Aristidou, and G.~Hug, ``Active
  distribution grids providing voltage support: The swiss case,'' \emph{IEEE
  Transactions on Smart Grid}, vol.~12, no.~1, pp. 268--278, 2021.

\bibitem{singhal_framework_2018}
A.~{Singhal} and V.~{Ajjarapu}, ``A framework to utilize ders’ var resources
  to support the grid in an integrated t-d system,'' in \emph{2018 IEEE Power
  Energy Society General Meeting (PESGM)}, Aug 2018, pp. 1--5.

\bibitem{arnold_optimal_2016}
D.~B. Arnold \emph{et~al.}, ``Optimal dispatch of reactive power for voltage
  regulation and balancing in unbalanced distribution systems,'' in
  \emph{{IEEE} {Power} and {Energy} {Society} {General} {Meeting}}, Jul. 2016,
  pp. 1--5.

\bibitem{farivar_equilibrium_2013}
M.~Farivar, L.~Chen, and S.~Low, ``Equilibrium and dynamics of local voltage
  control in distribution systems,'' in \emph{52nd {IEEE} {Conference} on
  {Decision} and {Control}}, Dec. 2013, pp. 4329--4334.

\bibitem{schweitzer_lossy_2020}
E.~Schweitzer, S.~Saha, A.~Scaglione, N.~G. Johnson, and D.~Arnold, ``Lossy
  distflow formulation for single and multiphase radial feeders,'' \emph{IEEE
  Transactions on Power Systems}, vol.~35, no.~3, pp. 1758--1768, 2020.

\bibitem{low_convex_relaxation}
L.~Gan and S.~H. Low, ``Convex relaxations and linear approximation for optimal
  power flow in multiphase radial networks,'' in \emph{2014 Power Systems
  Computation Conference}, 2014, pp. 1--9.

\bibitem{noauthor_ansi_2016}
\emph{{ANSI} {C}84.1: {American} {National} {Standard} for {Electric} {Power}
  {Systems} and {Equipment}—{Voltage} {Ratings} (60 {Hertz})}.\hskip 1em plus
  0.5em minus 0.4em\relax NEMA, Oct. 2016.

\bibitem{alok2019}
A.~K. Bharati and V.~Ajjarapu, ``Investigation of relevant distribution system
  representation with \uppercase{DG} for voltage stability margin assessment,''
  \emph{IEEE Transactions on Power Systems}, vol.~35, no.~3, pp. 2072--2081,
  2020.

\bibitem{sun_master_2015}
H.~Sun \emph{et~al.}, ``Master-{Slave}-{Splitting} {Based} {Distributed}
  {Global} {Power} {Flow} {Method} for {Integrated} {Transmission} and
  {Distribution} {Analysis},'' \emph{IEEE Trans on Smart Grid}, vol.~6, no.~3,
  pp. 1484--1492, May 2015.

\bibitem{alok_cosim_2021}
A.~K. Bharati and V.~Ajjarapu, ``Smtd co-simulation framework with helics for
  future-grid analysis and synthetic measurement-data generation,'' \emph{IEEE
  Transactions on Industry Applications}, pp. 1--1, 2021.

\end{thebibliography}

\end{document}